\documentclass[twoside,leqno]{article}

\usepackage[letterpaper]{geometry}

\usepackage{siamproceedings}

\usepackage[T1]{fontenc}
\usepackage{amsfonts}
\usepackage{graphicx}
\usepackage{epstopdf}
\usepackage{enumitem}
\usepackage{algorithmic}
\usepackage{amssymb}
\ifpdf
  \DeclareGraphicsExtensions{.eps,.pdf,.png,.jpg}
\else
  \DeclareGraphicsExtensions{.eps}
\fi

\newsiamremark{remark}{Remark}
\newsiamremark{hypothesis}{Hypothesis}
\crefname{hypothesis}{Hypothesis}{Hypotheses}
\newsiamthm{claim}{Claim}

\newtheorem{cond}[theorem]{Condition}

\usepackage{amsopn}

\usepackage[many]{tcolorbox}
\usepackage{stmaryrd}
\usepackage{pict2e,picture}
\usepackage{mathptmx}
\usepackage{eepic}
\usepackage{cprotect}

\usepackage{bbm}
\usepackage[many]{tcolorbox}
\usepackage{mathrsfs}
\usepackage{graphicx}

\usepackage{xcolor}
\definecolor{aquamarine}{rgb}{0.5, 1.0, 0.83}
\definecolor{aqua}{rgb}{0.0, 1.0, 1.0}
\usepackage{subcaption}
\usepackage{enumerate}

\usepackage{tikz}       % for tikz pictures
\usetikzlibrary{arrows}
\usetikzlibrary{chains}
\usetikzlibrary{mindmap}
\usetikzlibrary{backgrounds}
\usetikzlibrary{automata}
\usepackage{pgfplots}

\usepackage{todonotes}
\newcommand{\RvdH}[1]{\todo[color=magenta, inline]{Remco: #1}}

\DeclareSymbolFont{extraup}{U}{zavm}{m}{n}
\DeclareMathSymbol{\varheart}{\mathalpha}{extraup}{86}
\DeclareMathSymbol{\vardiamond}{\mathalpha}{extraup}{87}
\newcommand{\ensymboldefinition}{$\blacktriangleleft$}

\makeatletter
\newcommand{\fixed@sra}{$\vrule height 2\fontdimen22\textfont2 width 0pt\shortrightarrow$}
\newcommand{\shortarrow}[1]{%
  \mathrel{\text{\rotatebox[origin=c]{\numexpr#1*45}{\fixed@sra}}}
}
\makeatother

\newcommand{\remove}[1]{}
\newcommand{\ubar}[1]{\underline{#1}}
\newcommand{\ord}{\ensuremath{\mathrm{ord}}}

\newcommand{\eqn}[1]{\begin{equation} #1 \end{equation}}
\newcommand{\eqan}[1]{\begin{align} #1 \end{align}}
\newcommand{\sss}{\scriptscriptstyle}
\newcommand{\Var}{{\rm Var}}

\newcommand{\R}{{\mathbb R}}

\newcommand {\vep}{\varepsilon}

\newcommand{\indic}[1]{\mathbbm{1}_{\{#1\}}}

\newcommand{\nn}{\nonumber}

\newcommand{\be}{\begin{equation}}
\newcommand{\ee}{\end{equation}}

\newcommand{\vertex}{o}

\newcommand{\N}{\mathbb N}

\newcommand{\e}{{\mathrm e}}

\newcommand{\ra}{\rightarrow}

\newcommand{\prob}{\mathbb{P}}
\newcommand{\expec}{\mathbb{E}}

\newcommand{\Bin}{{\sf Bin}}

\newcommand{\FD}{F_{\sss D}}

\newcommand {\convd}{\stackrel{d}{\longrightarrow}}
\newcommand {\convp}{\stackrel{\sss \prob}{\longrightarrow}}
\newcommand {\convas}{\stackrel{a.s.}{\longrightarrow}}

\newcommand{\op}{o_{\sss \prob}}
\newcommand{\Op}{O_{\sss \prob}}

\newcommand{\Ecal}   {\mathcal{E}}

\newcommand{\conn}{\longleftrightarrow}

\newcommand{\Cmax}{\mathscr{C}_{\rm max}}
\newcommand{\Csec}{\mathscr{C}_{\rm sec}}
\newcommand{\cluster}{\mathscr{C}}

\newcommand{\Ver}{U}
\newcommand{\bfwit}{\boldsymbol{w}}
\newcommand{\bfdit}{\boldsymbol{d}}

\newcommand{\erdos}{Erd\H{o}s-R\'enyi random graph}

\newcommand{\GRGnw}{{\rm GRG}_n(\boldsymbol{w})}

\newcommand{\CLnw}{{\rm CL}_n(\bfwit)}

\newcommand{\NRnw}{{\rm NR}_n(\bfwit)}

\newcommand{\CMnd}{{\rm CM}_n(\boldsymbol{d})}
\newcommand{\CM}{{\rm CM}}
\newcommand{\SP}{{\rm SP}}
\newcommand{\CMpi}{\ensuremath{\mathrm{CM}_n(\boldsymbol{d},\pi_n(\lambda))}}

\newcommand{\CMndpi}{{\rm CM}_n(\boldsymbol{d}; \pi)}
\newcommand{\CMndlam}{{\rm CM}_n(\boldsymbol{d}; \pi_n(\lambda))}
\newcommand{\CMND}{{\rm CM}_N(\boldsymbol{d}(\pi))}

\newcommand{\CMnD}{{\rm CM}_n(\boldsymbol{D})}

\newcommand{\PA}{{\rm PA}}

\newcommand{\indeg}{{\rm in-deg}}

\newcommand{\myparagraph}[1]{\smallskip

\paragraph{\bf #1}}
\newcommand{\invisible}[1]{}

\newcommand{\seta}{\eta_s}

\newcommand{\perc}{\pi_n}
\newcommand{\percl}{\pi_n(\lambda)}
\newcommand{\percn}{\pi_n}
\newcommand{\cf}{c_{\sss \mathrm{F}}}
\newcommand{\dif}{\mathrm{d}}
\newcommand{\rNR}{\mathrm{NR}_n}
\newcommand{\bld}[1]{\boldsymbol{#1}}
\newcommand{\bw}{\bld{w}}
\newcommand{\sC}{\mathscr{C}}
\newcommand{\sCj}{\mathscr{C}_{\sss (j)}}
\newcommand{\sW}{\mathscr{W}}
\newcommand{\cGinf}{\mathscr{G}_{\sss \infty}}
\newcommand{\Unot}{\mathbb{U}^0_{\shortarrow{6}}}
\newcommand{\biS}{\mathbf{S}^\lambda}
\newcommand{\iS}{S^\lambda}
\newcommand{\refl}[1]{\ensuremath{\mathrm{refl}(#1)}}

\newcommand{\bN}{\mathbb{N}}
\newcommand{\bR}{\mathbb{R}}
\newcommand{\vC}{\mathbf{C}}
\newcommand{\clusterold}{\mathscr{C}}
\newcommand{\susceptibilitypi}[2]{s^\pi_{#2}(#1)}
\newcommand{\csusceptibilitypi}[2]{S^\pi_{#2}(#1)}

\newcommand{\clusterpi}[2]{\mathscr{C}^\pi_{#1}(#2)}

\newcommand{\set}[1]{\left\{#1\right\}}
\newcommand{\cG}{\mathscr{G}}

\begin{document}

\newcommand\relatedversion{}
\renewcommand\relatedversion{\thanks{The full version of the paper can be accessed at \protect\url{https://arxiv.org/abs/0000.00000}}} % Replace URL with link to full paper or comment out this line

\title{\Large Percolation on random graphs\relatedversion}
    \author{Remco van der Hofstad\thanks{Department of Mathematics and Computer Science, Eindhoven University of Technology (\email{r.w.v.d.hofstad@tue.nl}).}}

\date{}

\maketitle

% Copyright Statement
% When submitting your final paper to a SIAM proceedings, it is requested that you include
% the appropriate copyright in the footer of the paper.  The copyright added should be
% consistent with the copyright selected on the copyright form submitted with the paper.
% Please note that "20XX" should be changed to the year of the meeting.

% Default Copyright Statement
\fancyfoot[R]{\scriptsize{Copyright \textcopyright\ 20XX by SIAM\\
Unauthorized reproduction of this article is prohibited}}

% Depending on which copyright you agree to when you sign the copyright form, the copyright
% can be changed to one of the following after commenting out the default copyright statement
% above.

%\fancyfoot[R]{\scriptsize{Copyright \textcopyright\ 20XX\\
%Copyright for this paper is retained by authors}}

%\fancyfoot[R]{\scriptsize{Copyright \textcopyright\ 20XX\\
%Copyright retained by principal author's organization}}

%\pagenumbering{arabic}
%\setcounter{page}{1}%Leave this line commented out.

\begin{abstract}
Percolation is a model for random damage to a network. It is one of the simplest models that displays a {\em phase transition}: when the network is severely damaged, it falls apart in many small connected components, while if the damage is light, connectivity is hardly affected. We study the location and nature of the phase transition on random graphs. In particular, we focus on the connectivity structure close to, or below, criticality, where components display intricate scaling behaviour such that a typical connected component has bounded size, while the average and maximal connected component sizes grow like powers of the network size. 

We review the recent progress that has been made in two important settings: random graphs whose expected adjacency matrix is close to being {\em rank-1}, the most prominent examples being the configuration model and rank-1 inhomogeneous random graphs, and {\em dynamic random graphs}, i.e., random graphs that grow with time, such as uniform and preferential attachment models. Remarkably, these two settings behave rather differently. In all cases, the {\em inhomogeneity} of the underlying random graph on which we perform percolation is of crucial importance.
\end{abstract}

\section{Introduction.}
In percolation, edges and/or vertices are randomly removed from a connected graph $G=(V(G),E(G))$. Often, percolation is studied on infinite connected graphs, such as the hypercubic lattice or infinite trees. See the books by Grimmett \cite{Grim99}, Bollob\'as and Riordan \cite{BolRio06} and Kesten \cite{Kest82} for background. Percolation on infinite graphs is the simplest imaginable model that displays a {\em phase transition}. Suppose edges in the graph are removed independently with some probability $\pi\in[0,1]$.\footnote{Traditionally, $p$ is used as the edge-retention probability. Due to the potential conflicts with other usages of $p$ here, we choose for $\pi$ instead.} We think of the retention probability $\pi$ as the proportion of edges that escapes being damaged. When $\pi$ is small, the damage is substantial, and the graph falls apart into many small and finite components. In particular, while $G$ had an infinite connected component, after percolation, for $\pi$ small, no infinite connected component remains. On the other hand, when the damage is very little and $\pi$ is close to 1, the graph remains to have an infinite connected component even after the damage is inflicted. By monotonicity in the edge-retention probability, there is a {\em unique} $\pi_c$ above which infinite structures remain, while below it, no infinite connected component exists. Throughout the remainder of this paper, we will write `component' to mean `connected component'.

\myparagraph{Motivation for percolation on finite graphs.}
Real-world networks are finite, yet very large, graphs. Such networks can be affected by vertices and edges being down or attacked, rendering them out of service. This is sometimes called {\em attack} or {\em failure vulnerability}. Such a failure process could arise from vertices or edges randomly failing. Then, the key question is how much connectivity is preserved while vertices or edges are randomly knocked out. We will often assume that this is a truly random process, in the sense that the edges or vertices fail with equal probabilities independently of one another. Such an attack is sometimes called a {\em random attack}. 
%It could also be that the attacker has some knowledge about the network, and can deliberately knock out specific vertices, in a {\em deliberate attack}. Whether these settings are highly different is of practical interest, as it quantifies the effect of a knowledgeable opponent.

Real-world networks are generally {\em sparse} (see e.g.\ \cite{Hofs17, Hofs25, Newm10a}), meaning that their number of edges grows linearly with the network size. The relevant notion of convergence for sparse graphs is the highly powerful theory of \emph{local convergence} \cite{AldSte04,BenSch01,Hofs24}, which investigates the structure of finite graphs from the perspective of a {\em uniformly chosen vertex} called a {\em root}. The space of locally-finite rooted graphs, apart from isomorphisms, can be equipped with a local topology that turns it into a {\em Polish space}. Graph sequences converge locally when expectations of bounded and continuous functions of these rooted graphs converge to expectations with respect to a limiting probability distribution on rooted graphs. While this may seem a weak notion, it has proven to be very effective in analysing structural properties of large random graphs. For example, while the maximal component size is not described by the local limit, under one additional assumption, it is \cite{Hofs21b}, and we therefore refer to it as `almost local'.

\myparagraph{Percolation on finite graphs.}  
Let $G=(V(G),E(G))$ be a finite graph. For convenience, we denote the vertex set as $V(G)=[n]=\{1, \ldots, n\}$, while the edge set $E(G)$ can be general. Let $\pi$ be the edge-retention probability. For $u,v\in[n]$, we write that $u\conn v$ when there exists a path of edges that are retained in percolation, and let $\cluster_v^\pi=\{u\in[n]\colon u\conn v\}$ be the collection of  vertices that are connected to $v$ in percolation on $G$ with percolation parameter $\pi$. Percolation on finite graphs is inherently different from percolation on infinite transitive graphs in that a {\em unique critical value}, such as present on many infinite transitive graphs, cannot exist. Indeed, let $\prob_{\pi}$ denote the percolation measure on the finite graph $G$. Then, the dependence on the edge-retention probability $\pi$ is smooth, and $\prob_{\pi}(\Ecal)$ is a polynomial of finite degree for every event $\Ecal$. As a result, there cannot be a sharp phase transition, in the sense that if $\pi$ is close to critical and $\delta$ very small, then also $\pi+\delta$ will be close to critical. Further, on a finite graph, no infinite components can exist. This means that we have to redefine the notion of a critical value.
\smallskip

A {\em local} way of quantifying the critical behaviour is to take a {\em uniform} vertex $\Ver$, and to investigate the component $\cluster_{\Ver}^\pi$ of that vertex after percolating with percolation parameter $\pi$. When 
	$$
	\limsup_{k\rightarrow \infty}\limsup_{n\rightarrow \infty}\prob_{\pi_n}(|\cluster_{\Ver}^{\pi_n}|\geq k)=0,
	$$ 
the components in the percolated graph are small and we can think of $(\pi_n)_{n\geq 1}$ as being subcritical, while when 
	$$
	\liminf_{k\rightarrow \infty}\liminf_{n\rightarrow \infty}\prob_{\pi_n}(|\cluster_{\Ver}^{\pi_n}|\geq k)>0,
	$$ 
the sequence $(\pi_n)_{n\geq 1}$ is supercritical. Of course, the boundary between these two settings may not be sharp, and are generally referred to as {\em finite-size corrections}.

\myparagraph{Percolation on random graphs.}
In this paper, we focus on percolation on random graphs, which poses an additional difficulty due to the fact that the graphs are {\em random}, and thus we deal with {\em double randomness} due to percolation and the random graph on which it acts. As a result, we can think of the stochastic models on random graphs as processes in a disordered medium. Especially for percolation, this is somewhat paradoxical, since percolation often provides the random medium on which other processes, in particular random walks, act. For a more extended overview of percolation on random graphs, we refer to \cite[Part III]{Hofs25}.
%\smallskip

\invisible{Let us now continue by studying the local definition of a critical value in the random graph setting. The random graphs that we are interested in are locally tree-like in the large $n$ limit, and this limit does have a sharp phase transition. This allows us to identify the critical value, at least to first order, as the critical value of percolation on the random tree. The nice thing about this is that percolation on a branching process tree again gives rise to a branching process tree. In the context of random graphs that have a unimodular branching process tree as a local limit, as the configuration model and rank-1 random graphs (recall Theorems \ref{thm-LWC-CM-SP} and \ref{thm-LWC-GRG-SP}), it is not hard to see that the critical value of the limiting tree equals $\pi_c=1/\nu$, where $\nu$ is the expected value of the forward degree. Indeed, for this value, the expected offspring {\em after percolation} equals 1, and thus the unimodular branching process is critical (recall Theorem \ref{thm-survvsextncBP-SP}). This suggests that the critical behaviour of percolation on such random graphs occurs close to this value, and will constitute one of the main results in this chapter. We next discuss the critical behaviour in some more detail, by highlighting some critical exponents for critical percolation on random graphs.}

\myparagraph{The largest components in percolation on random graphs.}  
We let $\big(|\cluster_{\sss(i)}^\pi|\big)_{i\geq 1}$ denote the ordered collection of components, where we break ties (if they exist) arbitrarily. Thus, in particular, 
	\eqn{
	|\Cmax^\pi|=|\cluster_{\sss(1)}^\pi|=\max_{i\geq 1} |\cluster_{\sss(i)}^\pi|
	}
denotes the largest or maximal component size, and we write $|\Csec^\pi|=|\cluster_{\sss(2)}^\pi|$ for the second largest component size. Note that the vector $\big(|\cluster_{\sss(i)}^\pi|\big)_{i\geq 1}$ has a {\em random} number of coordinates, which equals the number of disjoint components, and it will be convenient to extend this finite vector to an infinite vector by appending infinitely many zeros at the end of it. 
\smallskip

In the above, we have defined critical behaviour in terms of the {\em local structure} of the component of a random vertex after percolation. This does not give us much information though about the largest component of the graph after percolation. The difficulty in dealing with the largest component is that it is an inherently {\em global} quantity, making it a difficult object to get our hands upon from a local perspective. We now discuss the percolation critical value from this global perspective.

We say that $(\pi_n)_{n\geq 1}$ (possibly depending on the network size $n$) is {\em asymptotically supercritical} when there exists $\vep>0$ independent of $n$ such that
	\eqn{
	\label{pi-asy-super}
	\liminf_{n\rightarrow \infty} \prob_{\pi_n}\big(|\Cmax^{\pi_n}|\geq \vep n\big)=1,
	}
and {\em asymptotically subcritical} otherwise. Thus, for an asymptotically supercritical sequence $(\pi_n)_{n\geq 1}$ the maximal component contains, with high probability, a {\em positive} proportion of the vertices, while $|\Cmax^{\pi_n}|/n\convp 0$ for a subcritical sequence $(\pi_n)_{n\geq 1}$. Here $\convp$ denotes convergence in probability.

\invisible{Clearly, this is a monotone definition in the sense that if $(\pi_n)_{n\geq 1}$ is an asymptotically supercritical sequence and $\pi'_n\geq \pi_n$ for every $n$, then also $(\pi_n')_{n\geq 1}$ is an asymptotically supercritical sequence). Similarly, if $(\pi_n)_{n\geq 1}$ is an asymptotically subcritical sequence and $\pi'_n\leq \pi_n$ for every $n$, then also $(\pi_n')_{n\geq 1}$ is an asymptotically subcritical sequence. This can be easily seen by a useful {\em coupling} of all percolation measures on the same probability space, which goes under the name of the {\em Harris coupling}\label{Harris-coupling}. For this, let $(U_e)_{e\in E(G)}$ denote a sequence of iid uniform random variables. Then we declare the edge $e$ $\pi$-retained when $U_e\leq \pi$. Clearly, the $\pi$-retained edges have distribution $\prob_{\pi}$, while if an edge is $\pi$-retained and $\pi'\geq \pi$, then the edge is also $\pi'$-retained (see Exercise \ref{exer-mon-pi}).}

\myparagraph{Critical value in large graph limits.}  In many cases, we see that 
	\eqn{
	\frac{|\Cmax^\pi|}{n}\convp \theta(\pi),
	}
for some $\theta(\pi)\in[0,1]$. Then, we can define the asymptotical percolation critical value also in terms of $\theta(\pi)$ as 
	\eqn{
	\label{pic-def-general}
	\pi_c=\inf\{\pi \colon  \theta(\pi)>0\}.
	}
When $|\Cmax^\pi|/n\convp \theta(\pi)>0$, the component $\Cmax^\pi$ is often {\em unique}, and then we called it the {\em giant}. In words, for $\pi<\pi_c$,  there is no giant in the percolated graph, while, for $\pi>\pi_c$, there exists a giant. In what follows, we will focus on finer properties of $|\Cmax^\pi|$ for $\pi$ close to the critical value $\pi_c$.

\myparagraph{Critical behaviour and finite-size effects in large graph limits.}  
While the notion of asymptotically supercritical sequences is useful to determine whether, after percolation, a giant component still exists, it is less useful to describe the critical behaviour for percolation on finite graphs. Indeed, for critical sequences $(\pi_n)_{n\geq 1}$, we would expect that the scaling of largest percolation components is quite intricate. A large body of work has recently appeared aimed at identifying the critical behaviour of percolation on various random and non-random finite graphs. Typically, an explicit sequence $(\pi_n)_{n\geq 1}$ is then chosen, and the behaviour of the largest percolation components determined. In each of these cases, the sequences involved are asymptotically subcritical. From this large body of work, we distill the following general picture of what it means to be barely subcritical, critical and supercritical: We say that a sequence $(\pi_n)_{n\geq 1}$ is {\em barely subcritical} when it is asymptotically subcritical, and there exists an $\vep>0$ such that
	\eqn{
	\label{pi-barely-sub}
	\liminf_{n\rightarrow \infty} \prob_{\pi_n}\big(|\Csec^{\pi_n}|/|\Cmax^{\pi_n}|\geq \vep \big)=1,
	}
while we say that a sequence $(\pi_n)_{n\geq 1}$ is {\em barely supercritical} when it is asymptotically subcritical and
	\eqn{
	\label{pi-barely-super}
	\frac{|\Csec^{\pi_n}|}{|\Cmax^{\pi_n}|}\convp 0.
	}
Finally, we say that $(\pi_n)_{n\geq 1}$ is {\em in the critical window} when it is asymptotically subcritical and there exist $0<a<b<\infty$, such that
	\eqan{
	\label{pi-crit-window}
	0<\liminf_{n\rightarrow \infty} \prob_{\pi_n}\big(|\Csec^{\pi_n}|/|\Cmax^{\pi_n}|\in[a,b]\big)
	&\leq \limsup_{n\rightarrow \infty} \prob_{\pi_n}\big(|\Csec^{\pi_n}|/|\Cmax^{\pi_n}|\in[a,b]\big)<1.
	}
A few remarks are in place here. In most cases, much stronger results are being proved in the literature. Often, in the barely subcritical regime, explicit sequences $(a_n)_{n\geq 1}$ are identified such that $|\cluster_{\sss(i)}^{\pi_n}|/a_n\convp c_i>0$ for every $i\geq 1$, which certainly implies \eqref{pi-barely-sub}. In the barely supercritical regime, often explicit sequences $(a_n)_{n\geq 1}$ are identified such that $|\Cmax^{\pi_n}|/a_n\convp c>0$, while $|\Csec^{\pi_n}|/a_n\convp 0$. Again this certainly implies \eqref{pi-barely-sub}. Finally, in the critical window, often results are proved that there exists a sequence $(a_n)_{n\geq 1}$ such that $a_n^{-1}\big(|\cluster_{\sss(i)}^{\pi_n}|\big)_{i\geq 1}$ converges in distribution to a certain (infinite) random vector in a suitable topology. We will see quite a few of such sharper results in what follows. 

\section{Percolation on rank-1 random graphs.}
\label{sec-rank-1-RGs}
We start this overview by discussing a class of models that we will call {\em rank-1 random graphs}. Note that each random graph $G_n=([n],E(G_n))$ of size $n$ can equivalently be defined as a collection of random variables $(X_{uv})_{u,v\in [n]}$, where $X_{uv}$ denotes the number of edges between vertices $u,v\in [n]$. We assume that we are in the {\em undirected} setting, for which $X_{uv}=X_{vu}$ for every $u,v\in [n]$. We allow for {\em multi-edges}, in that we allow $X_{uv}\geq 2$ so that $X_{uv}$ indicates the number of edges between $u,v\in [n]$. Further, we allow for $X_{vv}\geq 1$ for all $v\in[n]$, so that the random graph has self-loops. For convenience, we think of $X_{vv}$ as denoting {\em twice} the number of self-loops incident to $v\in[n]$. 
\smallskip

To explain the name rank-1 random graphs, we consider the {\em expected} adjacency matrix, given by $\expec[X_{uv}]$, where the expectation is over the randomness in the random graph, and we assume that there exists a sequence $(a_v)_{v\in[n]}$ such that 
	\eqn{
	\label{rank-1-RG-def}
	\expec[X_{uv}]\approx a_ua_v.
	}
Thus, the expected adjacency matrix is close to being of rank 1. We focus on two examples, the {\em configuration model} in Section \ref{sec-CM}, and {\em rank-1 inhomogeneous random graphs} in Section \ref{sec-rank-1-IRG}.

\subsection{Percolation on configuration models.}
\label{sec-CM}
In this section, we discuss percolation on configuration models. We start by introducing the model in Section \ref{sec-CM-model}, and discuss the percolation phase transition on it in Section \ref{sec-trans-CM}. In Sections \ref{sec-tau>3-CM}--\ref{sec-tau(2,3)-CM}, we describe the main results on the scaling limits of critical components.

\subsubsection{The configuration model.}
\label{sec-CM-model}
The configuration model is a model in which the degrees of vertices are fixed beforehand, and was introduced by Bollob\'as in \cite{Boll80b}; see \cite{Hofs17, Hofs24} for extensive overviews. Fix an integer $n$ that denotes the number of vertices in the random graph. Consider a sequence of degrees $\bfdit=(d_v)_{v\in[n]}$.\footnote{In many cases, the vertex degrees actually depend on $n$, and it would be more appropriate, but also more cumbersome, to write the weights as $\bfdit^{\sss(n)}=(d_v^{\sss(n)})_{v\in[n]}$. To keep notation simple, we refrain from making the dependence on $n$ explicit.} The aim is to construct an undirected multi-graph with $n$ vertices, where vertex $v$ has degree $d_v$. Without loss of generality, we assume that $d_v\geq 1$ for all $v\in [n]$, since  vertex $v$ is isolated when $d_v=0$ and
can thus be removed from the graph. We assume that the total degree
    	\eqn{
    	\ell_n=\sum_{v\in [n]} d_v
    	}
is even. To construct our multi-graph, we have $n$ separate vertices and, incident to vertex $v$, we have $d_v$ half-edges. Every half-edge needs to be connected to another half-edge to form an edge, and by forming all edges we build the graph. For this, the half-edges are numbered in an arbitrary order from $1$ to $\ell_n$. We start by randomly connecting the first half-edge with one of the $\ell_n-1$ remaining half-edges. Once paired, two half-edges form a single edge of the multi-graph, and the half-edges are removed from the list of unpaired half-edges. Hence, a half-edge can be seen as the left or the right half of an edge. We continue the procedure of randomly choosing and pairing the half-edges until all half-edges are connected, and call the resulting graph the {\it configuration model with degree sequence $\bfdit$}, abbreviated as $\CMnd$. 

Let us continue by explaining why we call the configuration model a rank-1 random graph. Since each half-edge incident to vertex $v$ is connected to a half-edge incident to vertex $u$ with probability $d_u/(\ell_n-1)$, we obtain
	\eqn{
	\expec[X_{uv}]=\frac{d_ud_v}{\ell_n-1},
	}
so that  \eqref{rank-1-RG-def} holds {\em with equality}, for $a_v=d_v/\sqrt{\ell_n-1}$. 

\myparagraph{Assumptions on the degrees.} To investigate the asymptotic properties of configuration models, we impose \emph{regularity conditions} on the degree sequence $\bfdit$. In order to state these assumptions, we introduce some notation. We denote the degree of a uniformly chosen vertex $\Ver$ in $[n]$ by $D_n=d_{\Ver}$. The random variable $D_n$ has distribution function $F_n$ given by
    \eqn{
    \label{def-Fn-CM}
    F_n(x)=\frac{1}{n} \sum_{j\in [n]} \indic{d_j\leq x},
    }
which is the {\em empirical degree distribution.}
We assume that the vertex degrees satisfy the following \emph{regularity conditions:}

\begin{cond}[Regularity conditions for vertex degrees]
\label{cond-degrees-regcond-SP}
~\\
{\bf (a) Weak convergence of vertex weights.}
There exists a distribution function $F$ such that, as $n\rightarrow \infty$,
    \eqn{
    \label{Dn-weak-conv}
    D_n\convd D,
    }
where $D_n$ and $D$ have distribution functions $F_n$ and $F$, respectively. Further, we assume that $F(0)=0$, i.e., $\prob(D\geq 1)=1$.\\
{\bf (b) Convergence of average vertex degrees.} As $n\rightarrow \infty$,
    \eqn{
    \label{conv-mom-Dn}
    \expec[D_n]\rightarrow \expec[D]\in (0,\infty),
    }
where $D_n$ and $D$ have distribution functions $F_n$ and $F$ from part (a), respectively.
%\\
%{\bf (c) Convergence of second moment vertex degrees.} As $n\rightarrow \infty$,
%    \eqn{
%    \label{conv-sec-mom-Dn}
%    \expec[D_n^2]\rightarrow \expec[D^2],
%    }
%where again $D_n$ and $D$ have distribution functions $F_n$ and $F$ from part (a), respectively, and the limit $\expec[D^2]$ may be infinite.
\end{cond}
Having defined the configuration model, we proceed by identifying its {\em percolation phase transition.}

\subsubsection{The phase transition of the configuration model.}
\label{sec-trans-CM}
In this section, we assume that Conditions \ref{cond-degrees-regcond-SP}(a)-(b) hold, and that $\CMnd$ is {\em supercritical}, in that 
 	\eqn{
	\label{sup-CM-perc}
   	\nu_n=\frac{\sum_{v\in [n]}d_v(d_v-1)}{\sum_{v\in [n]}d_v}\to \nu= \frac{\expec[D(D-1)]}{\expec[D]} >1.
  	}
In our first main result, which is due to Janson \cite{Jans09c}, we identify the percolation phase transition. In its statement, we write $\CMndpi$ for the result of percolation on $\CMnd$ with edge-retention probability $\pi$:

\begin{theorem}[Percolation phase transition in $\CMnd$ \cite{Jans09c}]
\label{thm-perc-PT-CM}
Suppose that Conditions \ref{cond-degrees-regcond-SP}(a)-(b) hold, and consider the random graph $\CMnd$, letting $n\to\infty$.
Assume that \eqref{sup-CM-perc} holds, so that $\CMnd$ is supercritical.
\begin{itemize}
\item[(a)]
For $\pi>1/\nu=\expec[D]/\expec[D(D-1)]$, there exists
$\theta(\pi)\in(0,1]$ such that $|\Cmax^\pi|/n\convp \theta(\pi)$,%=1-G_{\sss D}\big(1-\sqrt{\pi}+\eta^\star(\pi)\sqrt{\pi}\big),\\
%|E(\Cmax^\pi)|/n&\convp& \frac{1}{2}\expec[D] (1-\eta^\star(\pi)^2),
%\end{eqnarray*}
%where $\eta^\star(\pi)$ is the unique positive solution to 
%	\eqn{
%	\label{eta-star-def}
%	\eta^\star(\pi)=(1-\sqrt{\pi})+\sqrt{\pi}G_{D^\star-1}(1-\sqrt{\pi}+\eta^\star(\pi)\sqrt{\pi}).
%	}
while $|\Csec^\pi|/n \convp 0$.

\item[(b)]
For $\pi\leq 1/\nu$, $|\Cmax^\pi|/n\convp0$.
\end{itemize}
\end{theorem}

Theorem \ref{thm-perc-PT-CM} identifies the percolation critical value as $\pi_c=1/\nu=\expec[D]/\expec[D(D-1)]$, in the sense that any $\pi<\pi_c$ is asymptotically subcritical, while any $\pi>\pi_c$ is asymptotically supercritical. In particular, $\pi_c=0$ when $\expec[D^2]=\infty$, a property that is sometimes called {\em instantaneous percolation} or {\em robustness}. Robustness is desirable, since it implies that whichever the proportion of removed edges, there remains some decent connectivity in the percolated  graph.

We continue by investigating the {\em critical behaviour} of percolation on the configuration model.

\subsubsection{Scaling limits for finite third-moment degrees.}
\label{sec-tau>3-CM}
To formulate our main convergence result, for $p\geq 1$, define $\ell^p_{\shortarrow{6}}$ to be the set of infinite sequences 
$\boldsymbol{x}=(x_i)_{i\geq 1}$ with $x_1\geq x_2\geq \cdots\geq 0$ and $\sum_{i\geq 1} x_i^p<\infty$, and define the 
$\ell^2_{\shortarrow{6}}$ metric by 
    	\eqn{
    	d(\boldsymbol{x},\boldsymbol{y}) =\sqrt{\sum_{i\geq 1}(x_i-y_i)^2}.
    	}
Further, by $\ell^2_{\shortarrow{6}} \times \mathbbm{N}^{\infty}$, we denote the product topology of $\ell^2_{\shortarrow{6}}$ and $\mathbbm{N}^{\infty}$, where $\mathbbm{N}^{\infty}$ denotes the collection of sequences on $\mathbbm{N}$, endowed with the product topology. Define also
	\eqn{
	\mathbb{U}_{\shortarrow{6}}:= \big\{ ((x_i,y_i))_{i=1}^{\infty}\in  \ell^2_{\shortarrow{6}} \times \mathbbm{N}^{\infty}
	\colon \sum_{i=1}^{\infty} x_iy_i < \infty \text{ and } y_i=0 \text{ whenever } x_i=0 \; \forall i   \big\},
	} 
endowed with the metric 
	\eqn{
	\label{defn_U_metric}
	d_{\mathbb{U}}((\mathbf{x}_1, \mathbf{y}_1), (\mathbf{x}_2, \mathbf{y}_2)):= \bigg( \sum_{i=1}^{\infty} (x_{1i}-x_{2i})^2 \bigg)^{1/2}
	+ \sum_{i=1}^{\infty} \big| x_{1i} y_{1i} - x_{2i}y_{2i}\big|.
	} 
Finally, we introduce $\mathbb{U}^0_{\shortarrow{6}} \subset \mathbb{U}_{\shortarrow{6}}$ as 
	\eqn{
	\mathbb{U}^0_{\shortarrow{6}}:= \big\{((x_i,y_i))_{i=1}^{\infty}\in\mathbb{U}_{\shortarrow{6}} \colon \text{ if } x_k = x_m, k \leq m,\text{ then }y_k \geq y_m\big\}.
	}
For any $\mathbf{z}\in \mathbb{U}_{\shortarrow{6}}$, let $\ord(\mathbf{z})$ denote the element of $\mathbb{U}^0_{\shortarrow{6}}$ obtained by suitably ordering the coordinates of $\mathbf{z}$. The $\mathbb{U}^0_{\shortarrow{6}}$-topology will allow us to state joint convergence of component sizes and surpluses in a strong topology.
\smallskip

Consider bond percolation on $\CMnd$ with probability $\pi$, yielding  $\CMndpi$. Fix
  	\eqn{
	\label{pin-lambda-def}
  	\pi_{n}=\pi_n(\lambda):=\frac{1}{\nu_n} \bigg(1+ \frac{\lambda}{n^{1/3} }\bigg),
  	}
for some $\lambda\in \mathbbm{R}$.  Here we recall that $\nu_n=\expec[D_n(D_n-1)]/\expec[D_n]$, as defined in \eqref{sup-CM-perc}, and thus $\pi_n$ is close to $\pi_c=1/\nu$. The value $\pi_{n}$ as in \eqref{pin-lambda-def} rather than $\pi_c=1/\nu$ takes the finite-size corrections of the critical behaviour more closely into account. Note that $\pi_n(\lambda)$ in \eqref{pin-lambda-def} is non-negative for $n$ sufficiently large even when $\lambda<0$. 
\medskip

The main result in this section identifies the scaling limit of the ordered components of percolation on $\CMnd$ with edge-retention probability $\pi_n$ in \eqref{pin-lambda-def}. Let us now define the objects whose scaling limit we consider in detail. Recall that $(|\cluster_{\sss(j)}^\pi|)_{j\geq 1}$ denotes the ordered component sizes of $\CMndpi$. Let $\mathrm{SP}(\cluster)$ denote the surplus of the component $\cluster$. The main result of this section is then the following weak convergence result:

\begin{theorem}[Critical percolation on configuration models with finite third-moment degrees \cite{DhaHofLeeSen17}] 
\label{thm-percolation-CM}  
Consider $\CMnd$ with the degrees satisfying Condition \ref{cond-degrees-regcond-SP}(a)-(b) and \eqref{sup-CM-perc} with $\nu>1$, and that $\expec[D_n^3]\rightarrow \expec[D^3]<\infty.$ Fix $\pi_n=\pi_n(\lambda)$ as in \eqref{pin-lambda-def}. Then, there exists a non-degenerate limit $\mathbf{Z}(\lambda)$ such that, with respect to the $\mathbb{U}^0_{\shortarrow{6}}$ topology, as $n\rightarrow \infty$, 
    	\eqn{
    	\label{eqn_thm_percolation}
   	\mathbf{Z}_n(\lambda):=\ord\Big(\big(n^{-2/3}|\cluster_{\sss(j)}^{\pi_n}|, \mathrm{SP}(\cluster_{\sss(j)}^{\pi_n})\big)_{j\geq 1}\Big) \convd  \mathbf{Z}(\lambda).
  	}  
 \end{theorem}

\invisible{Next we consider the percolation component for multiple values of $\lambda$. There is a natural way to couple $(\CMndlam)_{\lambda\in \R}$ described as follows: for $\lambda<\lambda'$, perform bond-percolation on $\mathrm{CM}_n(\boldsymbol{d},\pi_n(\lambda'))$ with probability $\pi_n(\lambda)/\pi_n(\lambda')$. The resulting graph is distributed as  $\mathrm{CM}_n(\boldsymbol{d},\pi_n(\lambda))$. This can be used  to couple $(\mathrm{CM}_n(\boldsymbol{d},\pi_n(\lambda_i))_{i=0}^{k-1}$ for any fixed $k\geq 1$. The next theorem shows that the convergence of the component sizes holds jointly in finitely many locations within the critical window, under the above described coupling:

\begin{theorem}[Multiple times convergence]
\label{thm-multiple-convergence-CM}
Under the hypotheses of Theorem \ref{thm-percolation-CM} and with $\mathbf{C}_n(\lambda)=(n^{-2/3}|\mathscr{C}_{\sss (j)}(\lambda)| )_{j\geq 1}$, for any fixed $k\in \mathbbm{N}$ and $-\infty<\lambda_0< \lambda_1<\dots<\lambda_{k-1}<\infty$,
 	\eqn{
	\label{eqn_thm_multiple_convergence}
 	\big( \mathbf{C}_n(\lambda_0), \mathbf{C}_n(\lambda_1), \dots, \mathbf{C}_n(\lambda_{k-1})\big)  \convd \sqrt{\pi} (\boldsymbol{\gamma}^{ \lambda_0},	
	\boldsymbol{\gamma}^{ \lambda_1}, \dots, \boldsymbol{\gamma}^{\lambda_{k-1}})
 	}
with respect to the $(\ell^2_{\shortarrow{6}})^k$ topology where $\pi=1/\nu$. 
 \end{theorem}
\medskip}

\begin{remark}[Multiple values of $\lambda$, and informal relation to multiplicative coalescent]
\label{rem-mult-coal}
\normalfont The result in Theorem \ref{thm-percolation-CM} can be extended to {\em joint convergence} for several values of $\lambda$ at the same time. This process, with $\lambda\in \R$ considered as {\em time}, is a so-called {\em multiplicative coalescent process} \cite{Aldo97, AldLim98}, where components of masses $x_i(t)$ and $x_j(t)$ merge to form a component of mass $x_i(t)x_j(t)$ at rate $x_i(t)x_j(t)$. An intuitive picture is that as we change the value of the percolation parameter from $\pi_n(\lambda)$ to $\pi_n(\lambda+d\lambda)$, exactly one edge is added to the graph and the two endpoints $u,v$ are chosen approximately proportional to the number of half-edges of $u$ and $v$ that were not retained in $\CMndlam$. Define the \emph{degree deficiency} $\mathscr{D}_i$ of a component $\mathscr{C}_i$ to be the total number of half-edges in a component that were not retained in percolation. Think of $\mathscr{D}_i$ as the mass of $\mathscr{C}_i$. By the above heuristics, $\mathscr{C}_i$ and $\mathscr{C}_j$ merge at rate proportional to $\mathscr{D}_i\mathscr{D}_j$ and create a component of mass $\mathscr{D}_i+\mathscr{D}_j-2$. The proof with Dhara, van Leeuwaarden and Sen \cite{DhaHofLeeSen17} shows that the degree deficiency of a component is approximately proportional to the component size. Therefore, the component sizes merge \emph{approximately} like the multiplicative coalescent over the critical scaling window. \hfill \ensymboldefinition
\end{remark}

\subsubsection{Percolation on configuration models with infinite third-moment degrees.}
\label{sec-tau(3,4)-CM}
In the previous section, we have focused on the configuration model with finite third-moment degrees whose the percolation critical behaviour is in the same universality class as the \erdos{} (recall \cite{Aldo97}). This can be understood by the fact that, denoting the number of vertices found by pairing one half-edge by $d_{\sss(1)}-1$, we have that $\prob(d_{\sss(1)}=k)\approx k \prob(D_n=k)/\expec[D_n].$
Therefore, $\Var(d_{\sss(1)})\rightarrow \expec[D^3]/\expec[D]-(\expec[D^2]/\expec[D])^2<\infty$
when Condition \ref{cond-degrees-regcond-SP} holds and $\expec[D_n^3]\rightarrow \expec[D^3]$. Thus, the exploration process has {\em finite variance}, explaining why its scaling limit involves Brownian motion. In this section, we instead look at the heavy-tailed setting, where Condition \ref{cond-degrees-regcond-SP} holds, but $\expec[D^3]=\infty$, following work with Dhara, van Leeuwaarden and Sen \cite{DhaHofLeeSen17}. 

\myparagraph{Main results for heavy-tailed configuration models.}
Throughout this section, we use the shorthand notation
	\begin{subequations}
	\begin{equation}
	\label{eqn:notation-const}
	 \alpha= 1/(\tau-1),\qquad \rho=(\tau-2)/(\tau-1),\qquad \eta=(\tau-3)/(\tau-1),
	\end{equation}
	\begin{equation}\label{notation-power-tau(3,4)}
	a_n= n^{\alpha}L(n),\qquad b_n=n^{\rho}(L(n))^{-1},\qquad c_n=n^{\eta} (L(n))^{-2},
	\end{equation}
	\end{subequations}
where $\tau\in (3,4)$ and $L(\cdot)$ is a slowly-varying function. These sequences turn out to describe the sizes of the maximal degrees, the sizes of the largest components, and the width of the scaling window, respectively, when the degree distribution obeys a power-law with exponent $\tau\in(3,4)$ with slowly-varying corrections described by $n\mapsto L(n)$.
\smallskip

Before stating our main result, let us define the necessary conditions on the degrees under which it applies. Assume that there exists $\boldsymbol{\theta}=(\theta_1,\theta_2,\dots)\in \ell^3_{\shortarrow{6}}\setminus \ell^2_{\shortarrow{6}}$ such that
	\eqn{
	\label{def-degree-(3,4)}
 	\frac{d_i}{a_n}\to \theta_i \quad\text{ for all fixed }\quad i\geq 1,
	\qquad\text{and}\qquad 
	\lim_{K\to\infty}\limsup_{n\to\infty}a_n^{-3} \sum_{i=K+1}^{n} d_i^3=0,
	}
The first part of \eqref{def-degree-(3,4)} shows that $a_n$ in \eqref{notation-power-tau(3,4)} describes the order of the largest degrees, while the second part of \eqref{def-degree-(3,4)} shows that the largest degrees provide the main contribution to their third moment. The main result in this section describes the scaling of the largest critical percolation components on the configuration model with heavy-tailed degrees:

\begin{theorem}[Critical percolation on heavy-tailed configuration models \cite{DhaHofLeeSen20}]
\label{thm-perc-CM(3,4)}
Consider $\CMnd$ with the degrees satisfying Condition \ref{cond-degrees-regcond-SP}(a)-(b) and \eqref{sup-CM-perc} with $\nu>1$. Further, assume that \eqref{def-degree-(3,4)} holds, where $\boldsymbol{\theta}=(\theta_1,\theta_2,\dots)\in \ell^3_{\shortarrow{6}}\setminus \ell^2_{\shortarrow{6}}$. Assume that the percolation parameter $\pi_n=\pi_n(\lambda)$ satisfies
  	\begin{equation}
	\label{pic-CM-(3,4)}
  	\pi_n(\lambda):=\frac{1}{\nu_n} \big( 1+ \lambda c_n^{-1}+o(c_n^{-1}) \big)
  	\end{equation}
for some $\lambda\in \mathbbm{R}$. Then, there exists a non-degenerate limit $\mathbf{Z}(\lambda)$ such that, with respect to the $\mathbb{U}^0_{\shortarrow{6}}$ topology, as $n\to \infty$, 
	\begin{equation}
	\label{eqn:perc:limit}
 	\mathbf{Z}_n(\lambda) :=\ord \Big(\big( b_n^{-1}|\mathscr{C}_{\sss(j)}^{\pi_n}|,\mathrm{SP}(\mathscr{C}_{\sss(j)}^{\pi_n})\big)_{j\geq 1}\Big)\convd \mathbf{Z}(\lambda).
	\end{equation}
\end{theorem}

Note that $\tau\in(3,4)$, so that $b_n=n^{\rho+o(1)}$ with $\rho=(\tau-2)/(\tau-1)\in (\tfrac{1}{2},\tfrac{2}{3})$. Thus, the largest critical components of random graphs with {\em larger} degrees become {\em smaller}, a somewhat counterintuitive finding. It is related to the fact that the survival probability of a critical branching process with infinite-variance offspring is {\em smaller} than that of a branching process with finite-variance offspring. These connections are studied in quite some detail (under weaker conditions on the degrees) in work with Janson and Luczak \cite{HofJanLuc16}.

\myparagraph{Configuration model with iid degrees.} A natural setting to investigate the role of heavy-tailed degrees is to start with independent and identically distributed (iid) random degrees having a power-law distribution with $\tau\in(3,4)$. This is the setting investigated by Joseph \cite{Jose14}. We now relate this to the setting investigated so far. When dealing with iid degrees, it is helpful to {\em couple} the degrees in the configuration model for different values of $n$. Indeed, otherwise we cannot ever hope that the first convergence in \eqref{def-degree-(3,4)} holds for any fixed $i\geq 1$. For this, we use the fact that when $(E_i)_{i\geq 1}$ are iid exponentials, and $\Gamma_i=E_1+\cdots+E_i$ is the sum of the first $i$ of them, then the vector $(\Gamma_i/\Gamma_{n+1})_{i\in[n]}$ has the same distribution as the order statistics of $n$ iid uniforms $(U_i)_{i\in[n]}$ on $[0,1]$.
\smallskip

Let $\FD$ be the distribution function of the iid  degrees $(D_i)_{i\in[n]}$, so that $D_i$ has the same distribution as $[1-\FD]^{-1}(U_i)$. Let $(D_{\sss(i)})_{i\in[n]}$ be the order statistics, ordered such that $D_{\sss(1)}\leq D_{\sss(2)}\leq \cdots \leq D_{\sss(n)}$.
Then, we use that 
	\eqn{
	\label{coupling-degrees-tau(3,4)}
	\big(D_{\sss(n-i+1)}\big)_{i\in[n]} \stackrel{d}{=} \Big([1-\FD]^{-1}(U_{\sss(i)})\big)_{i\in[n]}  \stackrel{d}{=} \big([1-\FD]^{-1}(\Gamma_i/\Gamma_{n+1})\big)_{i\in[n]}.
	}
Equation \eqref{coupling-degrees-tau(3,4)} provides a highly powerful {\em coupling} between the degree sequences for random graphs of different sizes, which are such that in particular the {\em maximal degrees} are closely related. This mimics the setting in \eqref{def-degree-(3,4)}, and  the convergence result in it holds almost surely by the Strong Law of Large Numbers. This links the critical behaviour of $\CMnD$ to that of $\CMnd$ with deterministic degrees.

\subsubsection{Percolation on infinite-variance degree configuration models.}
\label{sec-tau(2,3)-CM}
Fix $\tau \in (2,3)$.  Recall \eqref{eqn:notation-const}. 
\invisible{For $v,t> 0$, define $M_t(v) := \sum_{j: v\theta_j\leq 1,\  t\theta_j\leq 1 } \theta_j^3.$
We will assume that for any $t>0$,
	\eqn{\label{density-assumption}
	\int_{0}^\infty  \e^{- tv^2M_t(v)} \dif v<\infty. 
	}
Also, let $D_n$ be the degree of a vertex chosen uniformly at random from $[n]$.} Assume that
	\begin{subequations}
	\begin{equation}
	\label{def-degree-(2,3)a}
 	n^{-\alpha}d_i\to \theta_i,
	\end{equation}
where $\boldsymbol{\theta}:=(\theta_i)_{i\geq 1}\in \ell^2_{\shortarrow{6}}\setminus \ell^1_{\shortarrow{6}}$ is such that $\theta_i = c_{\sss F}^\alpha i^{-\alpha}$, and 
	\eqn{
	\label{def-degree-(2,3)b}
	\lim_{K\to\infty}\limsup_{n\to\infty}n^{-2\alpha} \sum_{i=K+1}^{n} d_i^2=0.
	}
	\end{subequations}
We will show that the critical window for percolation on $\CM$ is given  by 
	\begin{equation}
	\label{eq:crit-window-CM}
	\pi_n=\pi_n(\lambda):= \frac{\lambda}{\nu_n}, \qquad \lambda \in (0,\infty),
	\end{equation}  
with $\nu_n$ as in \eqref{sup-CM-perc}. Notice that $\pi_n \sim n^{-2\alpha+1} \sim  n^{-\eta}$, where $\eta = (3-\tau)/(\tau-1)>0$ when \eqref{def-degree-(2,3)a}--\eqref{def-degree-(2,3)b} hold.
%The case where $\pi\ll \pi_c$ will be called the barely subcritical regime and the case $\pi_c\ll \pi \ll 1$ will be called the barely supercritical regime. 
%We will show that a {\em unique large component} emerges in the barely supercritical regime.
%We first state the results about the component sizes and the complexity in the critical window, and then discuss the barely sub-/supercritical regimes.

%We will always write $\sC_{\sss (i)}(\pi)$ to denote the $i^{\rm th}$ largest component in the percolated graph. The random graph on which percolation acts will always be clear from the context. 
A vertex is called {\em isolated} if it has degree zero in the graph $\CMpi$.  In the following theorem, which gives the asymptotics for the critical component sizes and the surplus edges of  $\CMpi$, it will be convenient to define the component size corresponding to an isolated vertex to be zero:
%For any component $\mathscr{C}\subset\CMpi$, let $\SP(\mathscr{C})$ denote the number of surplus edges given by $\# \{\text{edges in }\mathscr{C} \}-|\mathscr{C}|+1$. 

\begin{theorem}[Critical component sizes and surplus edges \cite{DhaHofLee21}]
\label{thm:main-CM(2,3)} 
Consider $\CMnd$ with the degrees satisfying Condition \ref{cond-degrees-regcond-SP}(a)-(b), and fix $\tau\in(2,3)$. Suppose that \eqref{def-degree-(2,3)a}--\eqref{def-degree-(2,3)b} hold. Fix $\pi_n=\pi_n(\lambda)$ as in \eqref{eq:crit-window-CM}. Then, there exists a non-degenerate limit  $\mathbf{Z}(\lambda)$, such that, with respect to the $\Unot$ topology, as $n\to\infty$, 
	\begin{equation}
	\mathbf{Z}_n(\lambda) := \ord\Big(\big(n^{-\rho}|\sCj^{\pi_n}|,\SP(\sCj^{\pi_n})\big)_{j\geq 1}\Big) \convd \mathbf{Z}(\lambda).
	\end{equation}
\end{theorem}

\invisible{\begin{remark}[Ignoring isolated components]
\label{rem:isolated} 
\normalfont
Note that $2\rho<1$ for $\tau\in (2,3)$. 
When percolation is performed with probability $p_c$, there are of the order $n$ isolated vertices and thus $n^{-2\rho}$ times the number of isolated vertices tends to infinity. 
This is the reason why we must ignore the contributions due to isolated vertices, when considering the convergence of the component sizes in the $\ell^2_{\shortarrow{6}}$-topology.
Note that an isolated vertex with self-loops does not create an isolated component. 
%Therefore, if we explore half-edges and .
\end{remark}}

\subsection{Percolation on rank-1 inhomogeneous random graphs.}
\label{sec-rank-1-IRG}
Having defined the main results for the configuration model, we move on to study rank-1 inhomogeneous random graphs, where the behaviour is remarkably similar.
We start by introducing the model in Section \ref{sec-IRG-model}, and then describe how results need to be changed in the finite-variance degree cases in Section \ref{sec-rank-1-fin-variance}. We close this section by discussing the surprising results for the case of rank-1 random graphs with infinite-variance degrees and single-edges in Section \ref{sec-NR-scale-free-single}.

\subsubsection{Rank-1 inhomogeneous random graphs.}
\label{sec-IRG-model}
In the \emph{Poissonian random graph} or Norros-Reittu model \cite{NorRei06}, the edge probability of the edge between vertices $u$ and $v$, for $u\neq v$, is equal to
    \eqn{
    \label{pij-NR}
    p_{uv} =p_{uv}^{\sss\rm(NR)}= 1-{\rm exp}\left(-w_u w_v/\ell_n\right),
    }
where $\bfwit=(w_v)_{v\in [n]}$ are the {\it vertex weights,} and now $\ell_n$ is the total weight of all vertices given by  $\ell_n=\sum_{v\in [n]} w_v.$
We denote the resulting graph by $\NRnw$.\footnote{In many cases, the vertex weights actually depend on $n$, and it would again be more appropriate to write the weights as $\bfwit^{\sss(n)}=(w_v^{\sss(n)})_{v\in[n]}$, and we again refrain from making the dependence on $n$ explicit.} The Poissonian random graph can be called rank-1, since $p_{uv}\approx w_u w_v/\ell_n$ when $w_u w_v/\ell_n$ is quite small. The  Poissonian random graph $\NRnw$ is close to many other inhomogeneous random graph models, such as the \emph{random graph with given prescribed degrees} or Chung-Lu model $\CLnw$ \cite{ChuLu02b,ChuLu06c, ChuLu06}, where instead
$p_{uv}=p_{uv}^{\sss\rm(CL)}=\min(w_uw_v/\ell_n, 1).$ A further adaptation is the so-called {\em generalised random graph model} $\GRGnw$ \cite{BriDeiMar-Lof05}, in which $p_{uv}=p_{uv}^{\sss\rm(GRG)}=\frac{w_uw_v}{\ell_n+w_uw_v}.$ See Janson \cite{Jans08a} or \cite[Chapter 6]{Hofs17} for conditions under which these random graphs are \emph{asymptotically equivalent}, meaning that all events have equal asymptotic probabilities. 
\smallskip

For $\NRnw$, the vertex weight $w_v$, for $v\in[n]$, is closely related to the {\em expected degree} of vertex $v$. This can be seen by noting that, with $d_v$ denoting the degree of vertex $v$, and $X_{uv}$ the indicator that the edge $\{u,v\}$ is present in $\NRnw$,
	\eqn{
	\expec[d_v]=\sum_{u\in [n]\colon u\neq v} \expec[X_{uv}]\approx \sum_{u\in [n]} p_{uv}\approx \sum_{u\in [n]} \frac{w_u w_v}{\ell_n}= w_v,
	}
when approximating the sum over $u\in[n]$ such that $u\neq v$ by the sum over all $u\in [n]$. While this is not exact, for large $n$ and most $v$, this is a good approximation. Thus, we are led to the conclusion that $\expec[d_v]$ is close to $w_v$, which gives a natural interpretation for the vertex weights. 

An alternative way of constructing $\NRnw$ is to independently consider $X_{uv}' \sim \mathrm{Poisson} (w_uw_v/\ell_n)$ edges between any pair of vertices $u,v\in[n]$, and then take $X_{uv}=\indic{X_{uv}'\geq 1}$ for all $u\neq v$, and $X_{vv}=0$ for all $v\in[n]$. This has the advantage that $\sum_{u\in[n]} X_{uv}' \sim \mathrm{Poisson} (w_v)$, and $\NRnw$ is obtained from the multi-graph with $(X_{uv}')_{1\leq u\leq v\leq n}$ by {\em collapsing} multi-edges and {\em removing} self-loops.

\subsubsection{Rank-1 inhomogeneous random graphs with finite-variance degrees.}
\label{sec-rank-1-fin-variance}
Similarly to the configuration model setting in Condition \ref{cond-degrees-regcond-SP}, we let $W_n=w_{\Ver}$ be the weight of a random vertex, and in our main results, we will make assumptions on $W_n$ similar to those for the configuration model:

\begin{theorem}[Scaling limit for finite-variance rank-1 inhomogeneous random graphs \cite{BhaHofLee09a,BhaHofLee09b, Turo09}]
\label{thm-perc-NR}
Consider the single- or multi-edge $\NRnw$. Assume that $W_n\convd W$ for some random variable $W$ satisfying $\nu=\expec[W^2]/\expec[W]>1$.
\begin{itemize}
\item[(a)] 
Fix $\pi_n=\pi_n(\lambda)$ as in \eqref {pin-lambda-def}, and assume further that $\expec[W_n^3]\rightarrow \expec[W^3]<\infty$. Then the scaling limit result in Theorem \ref{thm-percolation-CM} holds.

\item[(b)] Fix $\pi_n=\pi_n(\lambda)$ as in \eqref{pic-CM-(3,4)}, with $L(n)\equiv 1$. Assume further that $w_v=\cf(n/v)^{1/(\tau-1)}$ for some $\cf>0$. Then the scaling limit result in Theorem \ref{thm-perc-CM(3,4)} holds.
\end{itemize}
\end{theorem}
To appreciate Theorem \ref{thm-perc-NR}(b), suppose we take $\tau>2$, and consider the distribution function~$F$ satisfying $[1-F](w) = (\cf/w)^{\tau-1}$ for some $\cf>0$ and all $w>\cf$. Then, $w_v= [1-F]^{-1}(v/n)$ satisfies $w_v=\cf(n/v)^{1/(\tau-1)}$. Obviously, $F$ obeys a power law with exponent $\tau>2$, so that also $W$ does. In summary, Theorem \ref{thm-perc-NR} implies that $\NRnw$ is in the same universality class as $\CMnd$. Informally, we should think about the asymptotic degree distribution $D$ for $\CMnd$ being replaced by ${\rm Poisson}(W)$, where $W$ is the asymptotic weight distribution.

\invisible{\subsubsection{Scale-free rank-1 inhomogeneous random graphs: multi-edge setting.}
\label{sec-NR-scale-free-multi}
We next move on to the scale-free setting, for which $\tau\in(2,3)$. The next theorem identifies the barely supercritical behaviour, a setting not studied in 

\begin{theorem}[Barely supercritical multi-edge setting \cite{DhaHof24}]
\label{thm:asymp-multi}
Fix $\tau\in(2,3)$, and let $w_v=(n/\cf v)^{1/(\tau-1)}$ and $\percn =\lambda_n n^{-(3-\tau)/(\tau-1)}$ , where $\lambda_n \to\infty$ with $\lambda_n = o(n^{(3-\tau)/(\tau-1)})$. Then, as $n\to\infty$, 
	\eqn{
	\frac{|\sC_{\sss (1)}(\percn)|}{n\percn^{1/(3-\tau)} } \convp \mu\kappa^{1/(3-\tau)}, %\qquad  \frac{W_{\sss (1)}(\percn) }{n\percn^{1/(3-\tau)} } \pto \mu\kappa^{1/(3-\tau)},
	} 
where $\kappa = \cf^{\tau-2} \Gamma(3-\tau)$, and $\Gamma(\cdot)$ denotes the gamma function. Further, for any $j\geq 2$,   $|\sC_{\sss (j)}(\percn)| = \op(n\percn^{1/(3-\tau)})$. %and $W_{\sss (j)}(\percn) = \oP(n\percn^{1/(3-\tau)}) $. 
\end{theorem}

\begin{remark}[Critical scaling for $\rNR(\bw,\percn)$] \normalfont 
Since $\rNR(\bw)$ lies in the same universality class as the configuration model, one would expect that the components of $\rNR(\bw,\percn)$ exhibit critical behaviour with identical scaling limits and scaling exponents as in Theorem \ref{thm:main-CM(2,3)}, proved in \cite{DhaHofLee21}, when $\lambda$ is fixed. However, \cite{DhaHof24} does not prove this.
\hfill \ensymboldefinition
\end{remark}}

\subsubsection{Scale-free rank-1 inhomogeneous random graphs: single-edge setting.}
\label{sec-NR-scale-free-single}
We restrict here to the single-edge setting, the multi-edge setting should be closely related to that for the configuration model as studied in Theorem \ref{thm:main-CM(2,3)}. While not all results have been proved for the multi-edge $\NRnw$, the barely supercritical regime is identified in \cite{DhaHofLee21}. 

Let $\seta=(3-\tau)/2>\eta=(3-\tau)/(\tau-1)$, and define
	\begin{equation}
	\label{eq:scaling-window}
	\perc =\percl:= \lambda n^{-\seta}, 
	\qquad \text{for }\lambda \in (0,\lambda_c),
	\end{equation}
where $\lambda_c$ is given by 
	\begin{gather}
    	\lambda_c := \frac{\cf^{-1/\alpha}}{2}\sqrt{\frac{(3-\tau)\mu^{1/\alpha}}{A_\alpha}}, \qquad\quad \text{with}\qquad\quad
     	A_\alpha:= \int_0^\infty \frac{1-\e^{-z}}{z^{1/\alpha}} \dif z. \label{eqn:A-B-alpha-def}
	\end{gather}
Since $\seta>\eta$, $\perc$ in the single-edge case is {\em much larger} than $\perc$ in the multi-edge case. There turn out to be more surprising differences between the two cases, as we now describe in our two main results:

\begin{theorem}[Critical regime for single-edge $\mathrm{NR}_n(\bw, \percn)$ \cite{BhaDhaHof25}] 
\label{thm:main-crit}
Fix $\tau\in(2,3)$, and let $w_v=\cf(n/v)^{1/(\tau-1)}$. Consider the single-edge $\rNR(\bw,\percn)$ with $\percn=\percl$ as in \eqref{eq:scaling-window}, and $\lambda\in (0,\lambda_c)$ for $\lambda_c$ as in \eqref{eqn:A-B-alpha-def}.  Then, there exists a non-degenerate limit $(\sW_{\sss (j)}^{\infty}(\lambda))_{j\geq 1}$ such that, with respect to the $\ell^2_{\shortarrow{6}}$-topology, as $n\to \infty$, 
	\eqn{
	\label{eq:com-weight-crit-thm}
	(n^{\alpha}\percl)^{-1} (|\sCj^{\perc}|)_{j\geq 1}  \convd(\sW_{\sss (j)}^{\infty}(\lambda))_{j\geq 1}.
	%, \qquad
	%\text{and} \qquad n^{-\alpha} (W_{\sss (i)} (\percl))_{i\geq 1} \convd(\sW_{\sss (i)}^{\infty}(\lambda))_{i\geq 1} 
	}
\invisible{Moreover, for each $i\geq 1$, 
	\eqn{
    	(n^{\alpha}\percl)^{-1} |\sCj(\percl)| - n^{-\alpha} (W_{\sss (i)} (\percl) \convp 0,
	}
so that the convergence in \eqref{eq:com-weight-crit-thm} holds jointly.}
\end{theorem}
%The non-degenerate scaling limit of the component sizes, as well as their weights, is the hallmark of critical behaviour.
%\smallskip

Let us next consider percolation with probability $\percn = \lambda n^{-\seta}$ for $\lambda>\lambda_c$, for which we will show that a `tiny' giant component of size $\Theta(\sqrt{n})$ 
appears in the graph, and the size of this giant component concentrates. Moreover, the giant component is {\em unique} in the sense that the second largest component is of a smaller order.  To describe the limiting size of the giant component, fix $a>0$, and define 
	\eqn{
	\label{defn:zeta-lambda}
	\zeta_a^\lambda := \lambda\int_0^a  \cf u^{-\alpha} \rho_a^\lambda(u)\dif u,
	}
where $\rho_a^\lambda: (0,a] \to [0,1]$ is the maximum solution to the fixed point equation 
	\eqn{
	\label{defn:surv-prob}
	\rho_a^\lambda(u)=1-\e^{-\lambda \int_0^a \kappa(u,v)\rho_a^\lambda(v)\dif v}, \quad\qquad \text{with }\quad\qquad  \kappa(u,v) := 1-\e^{-\cf^{2} (uv)^{-\alpha}/\mu}. 
	}
It turns out that $\zeta^{\lambda} = \lim_{a\to\infty} \zeta_a^\lambda$ exists, and that $\zeta^\lambda \in (0,\infty)$ whenever $\lambda>\lambda_c$. 
We now state our result for the emergence of the tiny giant component for $\lambda>\lambda_c$:

\begin{theorem}[$\sqrt{n}$-asymptotics of size and uniqueness tiny giant \cite{BhaDhaHof25}]
\label{thm:supcrit-bd}
Fix $\tau\in(2,3)$, and let $w_v=\cf(n/v)^{1/(\tau-1)}$. Consider the single-edge $\rNR(\bw,\percn)$ with $\percn=\lambda n^{-\seta}$ for some $\lambda>\lambda_c$. 
Then, as $n\to\infty$, 
	\eqn{
	n^{-1/2}|\Cmax^{\percn}|\convp \zeta^\lambda, \qquad \text{and} \qquad n^{-1/2}|\Csec^{\percn}|\convp 0,
	}
where $\zeta^{\lambda} = \lim_{a\to\infty} \zeta_a^\lambda$ with $\zeta_a^\lambda$ given by \eqref{defn:zeta-lambda}. 
%Further, $\{v\colon w_v\geq  n^{1/2+\delta}\}\subseteq \sC_{\sss (1)}(\percn)$ with high probability for every $\delta>0.$
\end{theorem}
For $\lambda>\lambda_c$, the `tiny giant' component emerges inside the subset of vertices with weight of order at least $\sqrt{n}$. The tiny giant component in the whole graph consists primarily of the 1-neighbourhood of this tiny giant, which gives rise to \eqref{defn:zeta-lambda}. 
%The intuition is discussed in more detail in Section~\ref{sec:proof-outline}.

\subsection{Proofs for rank-1 random graphs.}
\label{sec-proofs-rank-1}
In this section, we give an informal overview of the proofs of our results. Except for Theorems \ref{thm:main-crit} and \ref{thm:supcrit-bd}, we focus on the configuration model $\CMnd$, the proofs for $\NRnw$ are similar in spirit.

\myparagraph{The local limit, and the critical value.} It is well-known that the configuration model, assuming Condition \ref{cond-degrees-regcond-SP}, converges locally in probability to a branching process having root offspring $D$, and offspring distribution $D^\star-1$ for every other vertex. Here $D^\star$ is the size-biased distribution of $D$, given by 
	\eqn{
	\label{size-biased-degree}
	\prob(D^\star=k)=\frac{k\prob(D=k)}{\expec[D]}.
	}
This is a so-called {\em unimodular} branching process, and since only the root has a different offspring distribution, the percolation critical value is equal to $\pi_c=1/\expec[D^\star-1]=1/\nu.$ This suggests, as in Theorem \ref{thm-perc-PT-CM}, that there is a giant for $\pi>1/\nu$, while there is no giant when $\pi\leq 1/\nu$. We next explain how such results can be proved.

\myparagraph{The percolation phase transition: Janson's construction.} A beautiful proof of the percolation phase transition in Theorem \ref{thm-perc-PT-CM} relies on the realisation by Janson \cite{Jans09c} that percolation on $\CMnd$ can be described in terms of another configuration model. On the \erdos{} with edge probability $\lambda/n$, when we perform percolation with parameter $\pi$, we obtain another \erdos{}, now with parameter $\pi \lambda/n$. When the \erdos{} is supercritical in that $\lambda>1$, we thus obtain that $\pi_c=1/\lambda$. Janson's construction extends this observation to the configuration model, for which it is a little more involved.
\smallskip

Fix a degree distribution $\bfdit$. For $v\in [n]$, let $d_v(\pi)=\Bin(d_v, \sqrt{\pi})$, where $\big(\Bin(d_v, \sqrt{\pi})\big)_{v\in[n]}$ are independent binomial random variables.
Further, let 
	\eqn{
	N^+=N^+(\pi)=\sum_{v\in [n]} \big[d_v-\Bin(d_v, \sqrt{\pi})\big],
	}
and let $N=N(\pi)=n+N^+$. For $v\in [N]\setminus [n]$, set $d_v(\pi)=1$ and color these vertices {\em red}. This has the following interpretation: We split a vertex of degree $d_v$ into a single vertex of degree $d_v(\pi)=\Bin(d_v, \sqrt{\pi})$, and $d_v-d_v(\pi)$ red vertices of degree 1. This is called the {\em explosion} of vertex $v$. We can think of the red vertices as being {\em artificial}, and they will later need to be removed. After this red-vertex removal, an edge is retained when both of its half-edges are retained, which occurs with probability $\sqrt{\pi}\times \sqrt{\pi}=\pi$, as it should be. Since all edge-retention decisions are made independently, this is thus equivalent to percolation on $\CMnd$. Let $\CMND$ denote the configuration model with vertex set $[N]$ and degree distribution $\bfdit(\pi)=(d_v(\pi))_{v\in [N]}$. Then, Janson's construction states that percolation on the configuration model $\CMndpi$ has the same distribution as $\CMND$ followed by the removal of the $N-n=N^+$ red vertices with degree 1. Here, we refer to the removal of a vertex of degree 1 as the act of removal of the vertex from the vertex set, together with the removal of its (single) edge from the edge set.

Janson's construction couples the configuration model and percolation on it on the same probability space. Since the removal of the red degree-1 vertices, and the edges incident to them, does not change the fact whether there exists a giant or not, the phase transition is thus the same as for $\CMND$. Let $D_n(\pi)$ denote $d_{\Ver}(\pi)$, where $\Ver\in[N]$ is chosen uniformly at random. Then, note that $\sum_{v\in [N]}d_v(\pi)=\sum_{v\in [n]}d_v$ since no half-edge is removed, while $\sum_{v\in [N]}d_v(\pi)(d_v(\pi)-1)=\sum_{v\in [n]}d_v(\pi)(d_v(\pi)-1)$, where $d_v(\pi)\sim\Bin(\ell_n, \sqrt{\pi})$ are independent binomial variables. It comes as no surprise that $\sum_{v\in [n]}d_v(\pi)(d_v(\pi)-1)$ is highly concentrated around its mean $\pi \sum_{v\in [n]}d_v(d_v-1)$ (where we use that $\expec[\Bin(m,p)(\Bin(m,p)-1)]=p^2m(m-1)$). Thus, $\nu_n(\pi)\convp \pi \nu$, and a giant can be expected to exist precisely when $\pi \nu>1$. The beauty of Janson's construction is that this result can be read off from the classical giant component result for the configuration model; see \cite[Part III]{Hofs25} for a discussion and references.

\myparagraph{Cluster scaling for critical configuration models with finite second-moment degrees.}
The proofs of Theorems \ref{thm-percolation-CM} and \ref{thm-perc-CM(3,4)} follow by again using Janson's construction. These proofs, alike the proof of Theorem \ref{thm-perc-PT-CM}, consists of three key steps. The first key step is the statement and proof of the corresponding result for the configuration model {\em without} percolation under appropriate conditions on the degree sequence. The second key step is to show that the percolated degree sequence in Janson's construction satisfies the conditions posed in the theorem for the critical configuration model. The final and third step is to describe the effect of removing the red vertices. In the remainder of this section, we will discuss how to analyse the scaling limit of the largest components for critical configuration models, i.e., without percolation. For Theorem \ref{thm-percolation-CM}, this means that we assume that
	\eqn{
	\label{nun-crit-CM}
   	\nu_n:= \frac{\expec[D_n(D_n-1)]}{\expec[D_n]} =1+\lambda n^{-1/3}+o(n^{-1/3}),
  	}
while for Theorem \ref{thm-perc-CM(3,4)}, $n^{1/3}$ is replaced by $c_n$.

We continue by describing the scaling limits in Theorems \ref{thm-percolation-CM}, \ref{thm-perc-CM(3,4)} and \ref{thm:main-CM(2,3)}, which all have a remarkably similar form.

\myparagraph{The scaling limit for finite third-moment degrees.}
Let us continue by defining the necessary objects so as to be able to describe the scaling limit of $\mathbf{Z}_n(\lambda)$. Define the Brownian motion with negative parabolic drift by
	\eqn{
	\label{def:inhomogen:BM}
	S^{\lambda}_{\mu,\kappa}(t)=\frac{\sqrt{\kappa}}{\mu} B(t) +\lambda t-\frac{\kappa t^2}{2\mu^{3}},
	}
where $\mathbf{B}= \big(B(t)\big)_{t \geq 0}$ is a standard Brownian motion, and $\mu>0$, $\kappa>0$ and $\lambda\in \R$ are constants to be determined later on. Define the reflected version of $\mathbf{S}^{\lambda}_{\mu,\kappa}=\big(S^{\lambda}_{\mu,\kappa}(t)\big)_{t \geq 0}$ as
	\eqn{
	\label{def-reflected-BM}
	\refl{S^\lambda_{\mu,\kappa}}(t) = S^{\lambda}_{\mu,\kappa}(t) - \min_{0 \leq s \leq t} S^{\lambda}_{\mu,\kappa}(s).
	}
In terms of the above, define $\gamma_{j}^\lambda$ to be the ordered excursions of the inhomogeneous Brownian motion $\mathbf{S}^\lambda_{\mu,\kappa}$ with
	\eqn{
	\label{parameter-perc-CM}
	\mu=\expec[D(\pi_c)],\qquad \kappa=\expec[D(\pi_c)^3]\expec[D(\pi_c)]-\expec[D(\pi_c)^2]^2, \quad \beta=1/\expec[D(\pi_c)],
	}
where $\pi_c=1/\nu$ is the asymptotic critical value as defined in Theorem \ref{thm-perc-PT-CM}. Denote  the $j^{\rm th}$ largest component of $\CMndlam$ by  $\mathscr{C}_{\sss (j)}(\lambda)$. We will also identify the scaling limit of the {\em surplus edges} of the large critical components. For this, define the counting process of marks $\mathbf{N}^\lambda= (N^\lambda_{\mu,\kappa}(s) )_{s \geq 0}$ to be a process that has intensity $\beta \refl{S^\lambda_{\mu,\kappa}}(s)$ at time $s$ conditional on $(\refl{S^\lambda_{\mu,\kappa}}(u) )_{u \leq s}$, so that
	\eqn{
	\label{def-counting-process}
	N^\lambda_{\mu,\kappa}(s) - \int\limits_{0}^{s} \beta \refl{S^\lambda_{\mu,\kappa}}(u)du
	}
is a martingale (see Aldous \cite{Aldo97}). %Alternatively, we can define $\mathbf{N}^\lambda$ to be a mixed-Poisson process where $N^\lambda(s)$ is Poisson with parameter $\int_0^s \refl{S^\lambda_{\mu,\kappa}}(u)du$.
\smallskip

For an excursion $\gamma$ starting at time $l(\gamma)$ and ending at time $r(\gamma)$, let $|\gamma|=r(\gamma)-l(\gamma)$ be its length, and $N(\gamma)=N^\lambda_{\mu,\kappa}(\gamma)$ denote the number of marks in the interval $[l(\gamma),r(\gamma)]$. 
%Then,  $( (| \gamma_j^\lambda | , N(\gamma_j^\lambda)))_{ j \geq 1} $ can be ordered as an element of  $\mathbb{U}^0_{\shortarrow{6}}$ almost surely by Bhamidi, Budhiraja and Wang \cite[Theorem 4.1]{BhaBudWan14}. Denote this element of $\mathbb{U}^0_{\shortarrow{6}}$ by 
Then, the scaling limit in Theorem \ref{thm-percolation-CM} is given by $\mathbf{Z}(\lambda)= \ord\Big(\big( \big(\big| \gamma_j^\lambda \big| , N(\gamma_j^\lambda)\big)\big)_{ j \geq 1}\Big) $.

\myparagraph{The scaling limit for infinite third-moment degrees.} Consider a decreasing sequence $ \boldsymbol{\theta}=(\theta_1,\theta_2,\dots)\in \ell^3_{\shortarrow{6}}\setminus \ell^2_{\shortarrow{6}}$. Denote $\mathcal{I}_i(s):=\indic{\xi_i\leq s }$ where $\xi_i\sim \mathrm{Exp}\Big(\frac{\theta_i}{\mu\sqrt{\nu}}\Big)$ independently, and $\mathrm{Exp}(r)$ denotes the exponential distribution with rate $r$.  Let $\mu=\expec[D]$. Consider the process 
	\begin{equation}
	\label{def-limiting-process-(3,4)}
	S^\lambda(t) =  \sum_{i=1}^{\infty} \frac{\theta_i}{\sqrt{\nu}}\left(\mathcal{I}_i(t)- \frac{\theta_i}{\mu\sqrt{\nu}}t\right)+\lambda t,
	\end{equation}
%Let $\tilde{\mathbf{S}}_{\infty}^{\lambda}$ denote the process in \eqref{def-limiting-process-(3,4)} with $\theta_i$ replaced by $\theta_i/\sqrt{\nu}$.
for some $\lambda\in\mathbbm{R},$ and define the reflected version of $S^\lambda(t)$ as in \eqref{def-reflected-BM}. Processes of the form \eqref{def-limiting-process-(3,4)} were termed \emph{thinned} L\'evy processes in work with Bhamidi and van Leeuwaarden \cite{BhaHofLee09b}, since the summands are thinned versions of Poisson processes. Indeed, we can interpret $\mathcal{I}_i(t)=\indic{N_i(t)\geq 1}$, where $(N_i(t))_{i\geq 1}$ are independent Poisson processes of rate $\theta_i/(\mu\sqrt{\nu})$, so that
	\eqn{
	\label{Levy-tau(3,4)}
	L^\lambda(t) = \sum_{i=1}^{\infty} \frac{\theta_i}{\sqrt{\nu}}\left(N_i(t)- \frac{\theta_i}{\mu\sqrt{\nu}}t\right)+\lambda t
	}
is a L\'evy process. Obviously, $S^\lambda(t)\leq L^\lambda(t)$. In what follows, we will explain that the indicator processes $(\mathcal{I}_i(t))_{i\geq 1}$ describe whether the vertex with the $i^{\rm th}$ largest degree has been found by the exploration or not, and then $\theta_i\mathcal{I}_i(t)/\sqrt{\nu}$ describes the jump made by the exploration process when $i$ has been found. Clearly, in the graph context, any vertex can only been found {\em once} for the first time, which explains why the process $(\mathcal{I}_i(t))_{i\geq 1}$ is relevant rather than the Poisson process $(N_i(t))_{i\geq 1}$. Thus, the fact that we deal with a thinned L\'evy process is directly related to a {\em depletion-of-points effect} present in exploring random graphs.
\smallskip

Let us continue to describe the necessary notion for the scaling limits of component sizes for critical percolation on $\CMnd$ in the infinite third-moment degree setting. For any function $f\in \mathbb{D}[0,\infty)$,  define $\ubar{f}(x)=\inf_{y\leq x}f(y)$.  $\mathbb{D}_+[0,\infty)$ is the subset of $\mathbb{D}[0,\infty)$ consisting of functions with positive jumps only. Note that $\ubar{f}$ is continuous when $f\in \mathbb{D}_+[0,\infty)$. An \emph{excursion} of a function $f\in \mathbb{D}_+[0,T]$ is an interval $(l,r)$ such that 
	\eqan{
	\label{def:excursion}
	\min\{f(l-),f(l)\}&=\ubar{f}(l)=\ubar{f}(r)=\min\{f(r-),f(r)\} \qquad \text{and}\qquad f(x)>\ubar{f}(r),\ \forall x\in (l,r)\subset [0,T].
	}
Excursions of a function $f\in \mathbb{D}_+[0,\infty)$ are defined similarly. We use $\gamma$ to denote an excursion, as well as the length of the excursion $\gamma$, to simplify notation. 
Then, as in \eqref{def-counting-process}, we define the counting process $\mathbf{N}$ to be the Poisson process that has intensity $\refl{S^\lambda}(t)$ at time $t$ conditional on $(S^\lambda(u) )_{u \leq t}$.
%, i.e., the process such that
%	\begin{equation} 
%	\label{def-counting-process-(3,4)}
%	N(t) - \int\limits_{0}^{t}\refl{S^\lambda}(u)du
%	\end{equation} 
%is a martingale (recall \eqref{def-counting-process} where a similar process was defined in the finite-third moment setting). 
We let $N(\gamma)$ denote the number of marks in the interval $\gamma$.
\smallskip

\invisible{Finally, we define a Markov process $(\mathbf{Z}(s))_{s\in\R}$ on $\mathbbm{D}(\R,\mathbb{U}^0_{\shortarrow{6}})$, called the {\em augmented multiplicative coalescent} (AMC) process. Think of  a collection of particles in a system with $\mathbf{X}(s)$ describing their masses and $\mathbf{Y}(s)$ describing an additional attribute at time $s$. Let $K_1,K_2>0$ be constants. The evolution of the system takes place according to the following rule at time $s$:
	\begin{itemize}
	\item[$\rhd$] For $i\neq j$, at rate $K_1X_i(s)X_j(s)$,  the $i^{\rm th}$ and $j^{\rm th}$ components merge and create a new component of mass 
	$X_i(s)+X_j(s)$ and attribute $Y_i(s)+Y_j(s)$.  
	\item[$\rhd$] For any $i\geq 1$, at rate $K_2X_i^2(s)$, $Y_i(s)$ increases to $Y_i(s)+1$. 
	\end{itemize}
Of course, at each event time, the indices are re-organized to give a proper element of $\mathbb{U}^0_{\shortarrow{6}}$. This process was first introduced by Bhamidi, Budhiraja and Wang in \cite{BhaBudWan14} to study the joint behaviour of the component sizes and the surplus edges over the critical window for so-called bounded-size-rules processes. Broutin and Marckert \cite{BroMar16} introduced it for the Erd\H{o}s-R\'enyi random graph.}%Bhamidi, Budhiraja and Wang \cite{BhaBudWan14} extensively study the properties of the  standard version of AMC, i.e., the case  $K_1=1,K_2=\tfrac{1}{2}$ and showed in \cite[Theorem 3.1]{BhaBudWan14} that this is a (nearly) Feller process.
%\smallskip

In terms of these quantities, the scaling limit $\mathbf{Z}(\lambda)$ in Theorem \ref{thm-perc-CM(3,4)} is given by 
 $\mathbf{Z}(\lambda):=\ord (\nu^{1/2}\gamma_i(\lambda),N(\gamma_i(\lambda)))_{i\geq 1}$,
where $\gamma_i(\lambda)$ is the $i^{\rm th}$ largest excursion of $S_{\infty}^{\lambda}$.

\myparagraph{The scaling limit for infinite second-moment degrees.}
Consider a decreasing sequence $ \boldsymbol{\theta}=(\theta_1,\theta_2,\dots)\in \ell^2_{\shortarrow{6}}\setminus \ell^1_{\shortarrow{6}}$. Denote $\mathcal{I}_i(s):=\indic{\xi_i\leq s },$ where $\xi_i\sim \mathrm{Exp}(\theta_i/\mu)$ independently, and consider the process 
	\begin{equation}
	\label{defn::limiting::process}
	\iS(t) =  \frac{\lambda \mu}{\|\bld{\theta}\|_2^2} \sum_{i=1}^{\infty} \theta_i\mathcal{I}_i(t)-  t,
	\end{equation}
for some $\lambda, \mu >0$, and define the reflected version $\refl{\iS}(t)$ of $\iS(t)$ as in \eqref{def-reflected-BM}. It turns out that, for any $\lambda>0$, the excursion lengths of the process $\biS = (\iS(t))_{t\geq 0}$ can be ordered almost surely as an element of $\ell^2_{\shortarrow{6}}$. We denote this ordered vector of excursion lengths by $(\gamma_i(\lambda))_{i\geq 1}$. Finally, again as in \eqref{def-counting-process}, define the counting process $\mathbf{N}^\lambda = (N^\lambda(t))_{t\geq 0}$ to be the Poisson process that has intensity $(\lambda\mu^2)^{-1}\|\bld{\theta}\|_2^2 \times\refl{ \iS}(t)$ at time $t$, conditionally on $( \iS(u) )_{u \leq t}$. 
%Formally, $\mathbf{N}^\lambda$ is characterised as the counting process for which 
%	\begin{equation} 
%	\label{defn::counting-process}
%	N^\lambda(t) - \frac{\|\bld{\theta}\|_2^2}{\lambda\mu^2} \int\limits_{0}^{t}\refl{ \iS(u)}\dif u
%	\end{equation} 
%is a martingale.  
We again let $N_i(\lambda)$ denote the number of marks of $\mathbf{N}^\lambda$ in the $i^{\rm th}$ largest excursion of $\biS$.
Then the scaling limit of critical component sizes in the case where $\tau\in(2,3)$ in Theorem \ref{thm:main-CM(2,3)} is equal to $\mathbf{Z}(\lambda):= \ord((\gamma_i(\lambda),N_i(\lambda)))_{i\geq 1}$.

\myparagraph{Exploration of components.} The fact that the scaling limits in Theorems \ref{thm-percolation-CM}, \ref{thm-perc-CM(3,4)} and \ref{thm:main-CM(2,3)} all have a similar shape is due to the fact that they all arise through an analysis of {\em exploration processes} of the components in the configuration model. We focus our explanation on Theorems \ref{thm-percolation-CM} and \ref{thm-perc-CM(3,4)}. 

Our \emph{breadth-first} exploration algorithm for the components in $\CMnd$ carries along vertices that can be alive, active, exploring and killed, and half-edges that can be alive, active or killed. We sequentially explore the graph as follows:
\begin{itemize}
\item[(S0)] At stage $i=0$, all the vertices and the half-edges are \emph{alive} but none of them are \emph{active}. Also, there are no \emph{exploring} vertices. 
\item[(S1)]  At each stage $i$, if there is no active half-edge at stage $i$, choose a vertex $v$ proportional to its degree among the alive (not yet killed) vertices, and declare all its half-edges to be \emph{active} and declare $v$ to be \emph{exploring}. If there is an active vertex but no exploring vertex, then declare the \emph{smallest} vertex to be exploring.
\item[(S2)] At each stage $i$, take an active half-edge $e$ of an exploring vertex $v$ and pair it uniformly to another alive half-edge $f$. Kill $e,f$. If $f$ is incident to a vertex $v'$ that has not been discovered before, then declare all the half-edges incident to $v'$ active (if any), except $f$. 
If $\mathrm{degree}(v')=1$ (i.e. the only half-edge incident to $v'$ is $f$) then kill $v'$. Otherwise, declare $v'$ to be active and larger than all other vertices that are alive. After killing $e$, if $v$ does not have another active half-edge, then kill $v$ also.

\item[(S3)] Repeat from (S1) at stage $i+1$ if not all half-edges are already killed.
\end{itemize}
Define the {\em exploration process} by
   	\begin{equation}
	\label{defn:exploration:process}
    	S_n(0)=0,\quad
     	S_n(l)=S_n(l-1)+d_{(l)}J_l-2,
    	\end{equation} 
 where $J_l$ is the indicator that a new vertex is discovered at time $l$ and $d_{(l)}$ is the degree of the new vertex chosen at time $l$ when $J_l=1$.  Let $\mathscr{C}_{k}$ be the $k^{\rm th}$ component explored by the above exploration process and define $\tau_{k}=\inf \big\{ i\colon S_{n}(i)=-2k \big\}.$ Then  $\mathscr{C}_{k}$ is discovered between the times $\tau_{k-1}+1$ and $\tau_k$, and  $\tau_{k}-\tau_{k-1}-1$ gives the total number of edges in $\mathscr{C}_k$. Since our critical components have small surplus, also $|\mathscr{C}_k|$ is close to $\tau_{k}-\tau_{k-1}-1$. This explains how {\em component sizes} are related to {\em excursion lengths} of appropriate stochastic processes.

\myparagraph{Scaling limit of the exploration process for finite third-moment degrees.} When $\expec[D_n^3]\rightarrow \expec[D^3]<\infty$, as in Theorem \ref{thm-percolation-CM}, $\expec[d_{(l)}]\approx \expec[D^\star]=2$, so that the expected jump sizes in \eqref{defn:exploration:process} have almost zero mean. Further, as explained at the start of Section \ref{sec-tau(3,4)-CM}, for any fixed $l$, $\Var(d_{\sss(l)})\rightarrow \expec[D^3]/\expec[D]-(\expec[D^2]/\expec[D])^2<\infty.$ Thus, we are dealing with a stochastic process that has almost zero mean, and finite variance, so that it will not come as a surprise that the scaling limit involves {\em Brownian motion}. Let us now relate this in more detail to the process in \eqref{def:inhomogen:BM}. We claim that 
	\eqn{
	\label{scaling-limit-exploration-fin-third}
	\Big(n^{-1/3} S_n\big(\lceil t n^{2/3}\rceil\big)\Big)_{t\geq 0}\convd (S^{\lambda}_{\mu,\kappa}(t))_{t\geq 0}.
	}
The fact that time is rescaled by $n^{2/3}$ is indicative of the fact that component sizes are of order $n^{2/3}$, while the rescaling by $n^{1/3}=\sqrt{n^{2/3}}$ corresponds to diffusive rescaling. Results of Aldous \cite{Aldo97} prove that \eqref{scaling-limit-exploration-fin-third} implies the scaling limit of component sizes in Theorem \ref{thm-percolation-CM}, while the surplus follows by tracking the number of $l$ for which $J_l=0$, which correspond to edges that are connected to active half-edges, and thus count surplus edges.

Let us next explain how \eqref{scaling-limit-exploration-fin-third} arise. The linear $\lambda t$ term in \eqref{defn:exploration:process} arises since the mean of our process is {\em not quite zero}, but is by \eqref{nun-crit-CM} equal to $\nu_n-1=\lambda n^{-1/3}+o(n^{-1/3})$. Thus, in $\lceil t n^{2/3}\rceil$ steps, the total drift is close to $\lambda n^{-1/3} \lceil t n^{2/3}\rceil\approx \lambda t n^{1/3}$, which after rescaling by $n^{-1/3}$ as in \eqref{scaling-limit-exploration-fin-third} gives rise to the linear drift term $\lambda t$ in \eqref{defn:exploration:process}. 

We next explain the occurrence of the negative parabolic drift in \eqref{defn:exploration:process}, which is more delicate. While it is true that $\expec[d_{(l)}-2]\approx \lambda n^{-1/3}+o(n^{-1/3})$ for {\em small} $l$, we actually need to know what $\expec[d_{(l)}-2]$ is for $l$ of the order $n^{2/3}$, for which it turns out that $\expec[d_{(l)}-2]$ becomes smaller. This is because larger vertices are more likely to be chosen early, due to the fact that vertices are chosen according to their degree. In order to describe this effect, given a sequence of weights $\bfdit=(d_v)_{v\in[n]}$, the {\em size-biased reordering} of $[n]$ is the random vector of vertices $(v_i)_{i\in [n]}$ such that 
	\eqn{
	\prob(v_i=j\mid v_1,\ldots, v_{i-1})=\frac{d_j}{\sum_{\ell\not\in \{v_1,\ldots, v_{i-1}\}} d_{\ell}}.
	}
It turns out that the sequence $(d_{(l)}J_l)_{l\geq 1}$, removing the zeros for which $J_l=0$ (i.e., the degrees of vertices that are actually new) have {\em exactly} the same distribution as $(d_{v_l})_{l\geq 1}$. The key analytical result is then that
	\eqn{
	\expec[d_{(\lceil t n^{2/3}\rceil)}-2]\approx n^{-1/3} \Big(\lambda-\frac{t (\sigma_{3}-2\sigma_2)}{\mu^2}\Big),
	\quad
	\text{ so that } 
	\quad
	\sum_{l=1}^{\lceil t n^{2/3}\rceil} \expec[d_{(l)}-2]\approx n^{1/3} \Big(\lambda t -\frac{t^2 (\sigma_{3}-2\sigma_2)}{2\mu^2}\Big).
	}
This explains how the negative quadratic drift in \eqref{def:inhomogen:BM} arises. The scaling limit in \eqref{scaling-limit-exploration-fin-third} follows by applying a Martingale Central Limit Theorem.

\myparagraph{Scaling limit of the exploration process for infinite third-moment degrees.} 
We next extend this analysis to the case where $\tau\in(3,4)$ as in Theorem \ref{thm-perc-CM(3,4)}, for which we conveniently rewrite \eqref{defn:exploration:process}. Call a vertex \emph{discovered} if it is either active or killed. Let $\mathscr{V}_l$ denote the set of vertices discovered up to time $l$ and $\mathcal{I}_i^n(l):=\indic{i\in\mathscr{V}_l}$. Note that 
	\begin{equation}
	\label{exploration-process-(3,4)}
    	S_n(l)= \sum_{i\in [n]} d_i \mathcal{I}_i^n(l)-2l=\sum_{i\in [n]} d_i \left( \mathcal{I}_i^n(l)-\frac{d_i}{\ell_n}l\right)+\left( \nu_n-1\right)l.
	\end{equation} 
Recall from \eqref{scaling-limit-exploration-fin-third} that the exploration process runs on scale $n^{1/3}$ when $\expec[D^3]<\infty$. For $\tau\in(3,4)$, instead, the maximal degrees in \eqref{def-degree-(3,4)} are of order $a_n$, which, by \eqref{eqn:notation-const}--\eqref{notation-power-tau(3,4)}, is equal to $n^{1/(\tau-1)+o(1)}$, where $1/(\tau-1)\in (\tfrac{1}{3},\tfrac{1}{2})$. Thus, even the addition of a {\em single} high-degree vertex makes $S_n(l)$ jump by more than $n^{1/3}$. It can thus be expected that the high-degree vertices are precisely what makes the exploration process move up macroscopically, while there is a natural drift downward due to the subtracted terms in \eqref{exploration-process-(3,4)}. Comparing \eqref{exploration-process-(3,4)} to \eqref{def-limiting-process-(3,4)}, we can expect that the linear term $\lambda t$ in \eqref{def-limiting-process-(3,4)} arises from the $(\nu_n-1)l$ term in \eqref{exploration-process-(3,4)}. Also the sum over $i$ in \eqref{def-limiting-process-(3,4)} has a clear relation to the sum over $i$ in \eqref{exploration-process-(3,4)}. The key scaling limit result for $\tau\in(3,4)$ is
	\eqn{
	\label{scaling-limit-exploration-(3,4)}
	\Big(a_n^{-1} S_n\big(\lceil t b_n\rceil\big)\Big)_{t\geq 0}\convd (S^{\lambda}(t))_{t\geq 0},
	}
with $(S^{\lambda}(t))_{t\geq 0}$ as in \eqref{def-limiting-process-(3,4)}. Since the time scale of the exploration process corresponds to the component sizes, this further explains how the $b_n$ scaling of the component sizes arises. This can also be understood from \eqref{notation-power-tau(3,4)}, since the {\em probability} that the exploration finds a vertex of degree of order $a_n$ is roughly equal to $a_n/\ell_n=\Theta(L(n)n^{\alpha-1})$, so that the {\em time} for a vertex of degree of order $a_n$ to be found is $\ell_n/a_n=\Theta(L(n)n^{1-\alpha})=\Theta(b_n)$, where the last equality follows since $\rho=1-\alpha$ by \eqref{eqn:notation-const}. The proof of \eqref{scaling-limit-exploration-(3,4)} crucially relies on the fact that the process in \eqref{exploration-process-(3,4)} can be uniformly approximated by restricting the sum to a {\em finite} number of terms, and the fact that $\big(\mathcal{I}_i^n\big(\lceil t b_n\rceil\big)_{t\geq 0, i\geq 1}\convd \big(\mathcal{I}_i(t)\big)_{t\geq 0, i\geq 1}$ in the product topology.

\myparagraph{Cluster scaling for critical configuration models with infinite second-moment degrees.}
The proof for $\tau\in(2,3)$ in Theorem \ref{thm:main-CM(2,3)} in \cite{DhaHofLee21} proceeds in a similar way as above, with one major difference, since the limiting process turns out to have {\em bounded variation almost surely}. Indeed, 
for any $u<t$, 
	\eqn{
	\expec\big[|\iS(t) - \iS(u)|\big] \leq \frac{\lambda \mu}{\|\bld{\theta}\|_2^2} \sum_{i=1}^{\infty} \theta_i \e^{-\theta_i u}(1-\e^{-\theta_i (t-u)})+ |t-u| \leq (\lambda\mu+1)|t-u|.
	}
However, since $\sum_i\theta_i = \infty$, the process experiences infinitely many jumps in any bounded interval of time. This means that the excursion theory, developed by Aldous and Limic \cite{AldLim98}, no longer applies, and thus has to be developed afresh.

\myparagraph{Tiny giant in single-edge scale-free Poissonian random graph.} To define the limiting variables in Theorem~\ref{thm:main-crit}, we introduce an infinite weighted random graph on $\N$ which belongs to a general class of models studied by Durrett and Kesten in \cite{DurKes90}. For this, we let vertex $v\in \N$ have weight $\theta_v:=\cf v^{-\alpha}\mu^{-1}$. Consider the random multi-graph $\cGinf(\lambda)$ on $\N$ where vertices $u$ and $v$ are joined independently by Poisson$(\lambda_{uv})$ many edges with $\lambda_{uv}$ given by 
	\begin{equation}\label{defn:lambda-ij}
	\lambda_{uv}:=\lambda^2\int_0^\infty \Theta_u(x) \Theta_v(x) \dif x, \quad \text{where}\quad
	\Theta_v(x):= 1-\e^{-\cf\theta_v x^{-\alpha}}.
	\end{equation}
For $i \geq 1$, let $\sW_{\sss (i)}^{\infty}(\lambda)$ denote the $i^{\rm th}$ largest element of the set $\big\{\sum_{i\in \sC}\theta_i\colon \sC \text{ is a component}\big\}.$
The following result implies that the limiting object is well-defined for $\lambda \in (0, \lambda_c]$, and undergoes a phase transition at $\lambda=\lambda_c$:

\begin{proposition}[Phase transition for the limiting model \cite{DurKes90}]
\label{prop-limit-as-finite}
For $\lambda\leq \lambda_c$, $ (\sW_{\sss (i)}^{\infty}(\lambda))_{i\geq 1}$ is in $\ell^2_{\shortarrow{6}}$ almost surely, while for $\lambda > \lambda_c$, $\cGinf(\lambda)$ is connected almost surely, in particular $\sW_{\sss (1)}^{\infty}(\lambda)=\infty$ and $\sW_{\sss (2)}^{\infty}(\lambda)=0$ almost surely.
\end{proposition}

The crux to the proof of Theorem \ref{thm:main-crit} in \cite{BhaDhaHof25} is that asymptotically there are Poisson($\lambda_{uv}$) many {\em two-step} paths between `macro-hubs' $u$ and $v$, via intermediate `meso-scale' hubs of weight $\Theta(n^\rho)$ in $\rNR(\bw, \perc(\lambda))$, for $u,v$ fixed as $n\rightarrow \infty$. 
The integral in \eqref{defn:lambda-ij} can be understood as the limit of the summation over the intermediate vertices in the two-step connection probabilities from $u$ to $v$.
These two-step paths between hubs form the backbone of the largest components. By Proposition \ref{prop-limit-as-finite}, the connectivity structure of these two-step connections undergoes a phase transition, and the resulting components on $\N$ are finite for $\lambda\leq \lambda_c$. This explains the critical phase for $\lambda\in (0,\lambda_c]$ in Theorem \ref{thm:main-crit}. 
\smallskip

The proof for $\lambda>\lambda_c$ is more delicate and technical. We only sketch the highlights. Fix a parameter $a>0$, and let 
	\eqn{
	\label{Nna-def}
	N_n(a)=\lfloor a n^{(3-\tau)/2} \rfloor,
	\qquad 
	\text{and}
	\qquad 
	N_n=N_n(1).
	}
Observe from \eqref{eq:scaling-window} that $\percn \asymp \lambda/N_n \asymp a\lambda/N_n(a)$. Since $w_v=\cf(n/v)^{\alpha}$ with $\alpha=1/(\tau-1)$ as defined in \eqref{eqn:notation-const},  for $u>0$,
	\eqn{
	\label{eq:asymp-w-N}
	w_{\lceil N_n u \rceil}
	=\cf u^{-\alpha}\Big(\frac{n}{N_n}\Big)^{\alpha}
	=\cf u^{-\alpha}\big(n^{(\tau-1)/2}\big)^{\alpha} \asymp \sqrt{n} \cf u^{-\alpha},
	}
and thus  $[N_n(a)]$ consists of vertices with weight at least of order $\sqrt{n}a^{-\alpha}$.  The key conceptual step is that, for large enough $a$, a giant component emerges inside $[N_n(a)]$ that forms the core connectivity structure of the tiny giant component in the whole graph. In turn, this graph is an inhomogeneous random graph, for which the critical value can be determined exactly, as we explain in more detail now.

To this end, consider the percolated graph  $\rNR(\bw,\percn)$, restricted to $[N_n(a)]$, and denote this subgraph by $\cG_{\sss N_n(a)}$. Then, $\cG_{\sss N_n(a)}$ is distributed as an inhomogeneous random graph that is {\em sparse} in that the number of edges grows linearly in the number of vertices in the graph. Thus, the emergence of the giant component within $\cG_{\sss N_n(a)}$ can be studied using the general setting of inhomogeneous random graphs developed by Bollob\'as, Janson and Riordan in \cite{BolJanRio07}. 
In particular, the results of \cite{BolJanRio07} gives a critical value ${\lambda_c(a)}$, such that, for $\lambda>\lambda_c(a)$, a unique giant with concentrated size exists inside $[N_n(a)]$, that is {\em stable} to the addition of a small proportion of edges. The stability result is crucial to understand the perturbation on this giant after adding all the edges due to connections outside $[N_n(a)]$.

It turns out that $\lim_{a\to\infty}{\lambda_c(a)}  = \lambda_c$, where $\lambda_c$ is given by \eqref{eqn:A-B-alpha-def}. The connection between ${\lambda_c(a)}$ and $\lambda_c$ is remarkable given their vastly different descriptions, and is proved by an explicit computation. The convergence of ${\lambda_c(a)}$ is also a key conceptual step, since it shows that, whenever $\lambda>\lambda_c$, one can choose $a$ to be large enough to make a giant appear inside~$[N_n(a)]$. Finally, the asymptotics for the giant inside $\cG_{\sss N_n(a)}$ are given by properties of certain multi-type branching processes that depend sensitively on $a$. Since $a>0$ is arbitrary, we need to take the limit of $a\rightarrow \infty$, and show that the size of the tiny giant in $[n]$ converges. For this, we identify the primary contributions to the size of the giant in the whole graph using the giant inside $[N_n(a)]$, for $a$ large enough, which consist of the vertices at {\em finite} distance from $[N_n(a)]$ in the percolated graph. To explain that these are of order $\sqrt{n}$, we note that a vertex of weight $w_v\gg 1$ has degree roughly $w_v \pi_n$ in the percolated graph. Thus, the {\em neighbourhood} of the giant in $[N_n(a)]$ in the percolated graph $\rNR(\bw,\percn)$ has order
	\eqn{
	\sum_{v=1}^{N_n(a)} w_v \pi_n=\Theta(1) n^{\alpha} \sum_{v=1}^{N_n(a)} v^{-\alpha} n^{-(\tau-3)/2}=\Theta(1) n^{\alpha} N_n^{1-\alpha}  n^{-(\tau-3)/2}=\Theta(\sqrt{n}),
	}
by \eqref{eq:scaling-window} and since $N_n=\Theta(n^{(3-\tau)/2})$ by \eqref{Nna-def}. This explains why the size of the tiny giant is of order $\sqrt{n}$. Its concentration and uniqueness follow from the fact that the giant in $[N_n(a)]$ is concentrated and unique, as proved in \cite{BolJanRio07}.

\section{Dynamic random graphs.}
\label{sec-dynamic-RGs}
In this section, we discuss recent progress on percolation on various {\em dynamic} random  graphs models. Such models can be seen as random graph models that {\em grow} in time, so that the graph sequence $(G_n)_{n\geq 1}$ forms a {\em stochastic process}. This is not the case for, e.g., $\CMnd$, for which the graphs of size $n$ and $n+1$ are unrelated. It turns out that the percolation phase transition is remarkably different for such dynamic random graphs. 

We study preferential and uniform attachment random graphs in Section \ref{sec-UA}, and related dynamic random graphs (such as dynamic inhomogeneous random graphs) in Section \ref{sec-related-dynamic-random-graphs}.

\subsection{Percolation on preferential and uniform attachment random graphs.}
\label{sec-UA}
\subsubsection{Preferential and uniform attachment models.} Fix an integer parameter $m\geq 2$, weight parameters $a\in [0,1]$, $\delta>-m/a$, and consider an attachment function $f\colon \bN_0 \to \bR_+$,
\eqn{
    \label{eqn:linear-att-def}
    f(x) = a x + \delta, \qquad x\in \bN_0.
    }
Let $(G_n)_{n\geq 1}$ be a sequence of networks grown using preferential attachment dynamics and number of incoming edges $m$. Thus,  having constructed the network $G_n$ and writing $\deg(v,n)$ for the degree of an existing vertex $v$ in $G_n$,  a new vertex $n+1$ enters the system with $m$ edges which it uses to connect to the existing network, via sequentially connecting each edge to an existing vertex $v\in G_n$ in the network with probability approximately proportional to $f(\deg(v,n))$. In more detail, we start from any graph on $2$ vertices and finitely many connections between them such that at least one of the initial vertices has degree at most $m$. Let $d_1$ and $d_2$ denote the degrees of vertices $1$ and $2$ in the initial graph. Without loss of generality, consider $d_2\leq m$. For every new vertex $v$ joining the graph and $j=1,2,\ldots,m$, the attachment probabilities are given by
\eqn{\label{def:model:d}
	\prob\Big( v\overset{j}{\rightsquigarrow} u\mid\PA^{(d)}_{v,j-1}(m,\delta) \Big)= \frac{d_u(v,j-1)+\delta}{c_{v,j}^{(d)}}\hspace{1.5cm}\text{for}~u<v~,
}
where $v\overset{j}{\rightsquigarrow} u$ denotes that vertex $v$ connects to $u$ with its $j^{\rm th}$ edge, $\PA^{(d)}_{v,j}({m},\delta)$ denotes the graph on $v$ vertices, with the $v^{\rm th}$ vertex having $j$ out-edges, and $d_u(v,j)$ denotes the degree of vertex $u$ in $\PA^{(d)}_{v,j}({m},\delta)$. We identify $\PA^{(d)}_{v+1,0}({m},\delta)$ with $\PA^{(d)}_{v,m}({m},\delta)$.  The normalizing constant $c_{v,j}^{(d)}$ in \eqref{def:model:d} equals
	\eqn{
	\label{eq:normald}
	c_{v,j}^{(d)} := a d_{[2]}+2\delta+(2am+\delta)(v-3)+a(j-1),
	}
where $d_{[2]}=d_1+d_2$. We denote the above model as $\PA^{(d)}_{v}({m},\delta)$, which is equivalent to $\PA_v^{(m,\delta)}(d)$ as defined in \cite{Hofs17} when we set $a=1$. The case where $a=0$ corresponds to the {\em uniform attachment} case. The main quantity of interest is percolation on the constructed network. Fix a parameter $\pi \in (0,1)$, and let $G^{\pi}_n$ denote the network obtained from $G_n$ by independently retaining each edge with probability $\pi$ and deleting it with probability $1-\pi$.

%While the main model studied in this paper will be the \emph{uniform attachment model},  we describe the general model as one major goal of the paper is to develop general stochastic approximation tools required to understand percolation on general models of network evolution.

\subsubsection{The phase transition.}
We first describe the phase transition, by identifying for which percolation parameter $\pi\in[0,1]$ there exists a giant component, and for which there does not:
\begin{theorem}[Critical percolation threshold for $\PA$ models \cite{HazHofRay23}]
\label{thm-PT-PAM}
	%Let $\{G_n\}_{n\geq 1}$ be a sequence of preferential attachment model ((a), (b) or (d)) with parameters $m$ and $\delta$
	%and $\mu$ be the law of the local limit P\'olya point tree with parameters $m$ and $\delta$.
	Fix $a=1$ and $m\geq 2$. For any $\pi\in[0,1]$, let $\clusterpi{\max}{n}$ and $\clusterpi{\rm sec}{n}$ denote the maximal and second largest component of $G_n^\pi$.
	%and $\mathscr{C}(\pi)$ denote the component in the percolated P\'olya point tree containing the root 
	Then
	\eqn{\label{eq:main:theorem}
	\frac{|\clusterpi{\max}{n}|}{n} \convp\theta(\pi)
	%=\mu\Big( |\mathscr{C}(\pi)|=\infty \Big)
	\qquad\text{and}\qquad\frac{|\clusterpi{\rm sec}{n}|}{n} \overset{\prob}{\to}0.
	}
	%where $\overset{\prob}{\to}$ denotes convergence in probability with respect to both the random graph and percolation. 
	Furthermore, define $\pi_c=0$ for $\delta\in(-m,0]$, and
	\eqn{
	\label{def-pi-c-PA}
	\pi_c=\frac{\delta}{2\big(m(m+\delta)+\sqrt{m(m-1)(m+\delta)(m+1+\delta)}\big)}\qquad \text{for }\delta>0.
	}
	Then
	$\theta(\pi)>0$ for $\pi>\pi_c$, whereas $\theta(\pi)=0$ for $\pi\leq \pi_c$.
\end{theorem}
It turns out that $1/\pi_c$ denotes the exponential growth of the generations of the local limit of the preferential attachment model, which is the P\'olya point tree (PPT) \cite{BerBorChaSab14, GarHazHofRay22}. Thus, for $\pi<\pi_c$, the exponential growth is removed, and the percolated PPT dies out almost surely, while for $\pi>\pi_c$, it continues to grow exponentially with positive probability. This explains the relevance of $\pi_c$ in \eqref{def-pi-c-PA}. Bollob\'as and Riordan \cite{BolRio05}, by adapting argument for the {\em uniformly grown random graph} \cite{BolJanRio05}, have shown that, for the uniform attachment case for which $a=0$ and $\delta=1$, the critical percolation threshold becomes $\pi_c=1/[2\big(m+\sqrt{m(m-1)}\big)]$.

\subsubsection{Subcritical percolation on the uniform attachment model.}
Here we restrict to the uniform attachment case, for which $a=0$, and $m=2$. Write $\clusterpi{v}{n}$ for the component of vertex $v$ in $G^{\pi}_n$.  % Finally let $\cC(i,n)$ denote the size of the component containing vertex $i$ at time $n$. 
Define 
\begin{equation}
\label{pic-UM-m=2}
    \pi_c = \frac{2 - \sqrt{2}}{4} = \frac{1}{2(2 + \sqrt{2})}.
\end{equation}
%Known results and proof techniques for percolation on preferential attachment models, as proved e.g.\ in \cite{HazHofRay23}, imply that the size of the maximal component satisfies $|\clusterold_{\sss(1)}^\pi(n)|/n \convp \theta(\pi)$, where $\theta(\pi)$ is strictly positive limit for $\pi > \pi_c$ and zero for $\pi < \pi_c$. Thus, $\pi < \pi_c$ is the {\em subcritical regime}. 
Further, for $\pi < \pi_c$, define 
\begin{equation}
\label{eqn:alpha-pi}
    \alpha(\pi) := \frac{1}{2} \left( 1 - \sqrt{8\pi^2 - 8\pi + 1} \right).
\end{equation}
We now state the main results in \cite{BanBhaHofRay25}:

% Write $\clusterpi{\sss \geq i}{n}$ via 
% \begin{equation}
% \label{eqn:cc-less-def}
% \clusterpi{\sss \geq i}{n} =
% \left\{
% \begin{array}{ll}
% \clusterpi{i}{n} & \text{if $i$ is the oldest vertex in the component,}  \\
% \emptyset & \text{otherwise}.
% \end{array}
% \right.
% \end{equation}

\begin{theorem}[Convergence of maximal components \cite{BanBhaHofRay25}]
\label{thm-max-comp}
Fix $a=0$ and $m=2$. Let $\pi<\pi_c$. 
\begin{itemize}
\item[(a)] For any fixed $i\in \N$, $n^{-\alpha(\pi)}|\clusterpi{i}{\pi}| \convas \zeta_i$  as $n\to\infty$, where  $\zeta_i >0$ almost surely.
\item[(b)] $n^{-\alpha(\pi)} |\clusterpi{\max}{n}| \convas \max_{i \ge 1}\zeta_i$ as $n\to\infty$. Further, the maximal component is `weakly persistent' in the sense that, for any $\vep>0$, there exists $K(\vep)$ such that 
	\begin{equation*}
    	\prob\left(\clusterpi{\max}{n} = \clusterpi{i}{n} \text{ for some } i> K(\vep) \text{ for infinitely many } n\right) < \vep.
	\end{equation*}
\end{itemize}
\end{theorem}

\subsection{Percolation on related dynamic random graphs.}
\label{sec-related-dynamic-random-graphs}
In Section \ref{sec-UA}, we have discussed the phase transition on preferential and uniform attachment models, and the subcritical phase on uniform attachment models. We believe that similar results hold for various related dynamic random graph models. M\"orters and Schleicher \cite{MorSch25} investigate a class of dynamic inhomogeneous random graphs, in which an edge is independently present between vertices $u$ and $v$ with probability $p_{uv}=\beta (u\wedge v)^{\gamma-1}(u\vee v)^{\gamma}$ for some $\gamma\in [0,1]$. Such models are called {\em inhomogeneous random graphs of preferential attachment type}, since also for preferential attachment models, the probability that an edge is present between vertices $u$ and $v$ is of order  $p_{uv}$ for $\gamma =(m+\delta)/(2m+\delta)$ and an appropriate $\beta>0$ (see \cite[Chapter 8]{Hofs24}). By an analysis as in \cite{DerMor13}, it can be shown that $\beta_c=(\tfrac{1}{4}-\tfrac{\gamma}{2})\wedge 0$. Thus, $\beta_c>0$ when $\gamma<\tfrac{1}{2}$, which corresponds to $\delta>0$. Then, M\"orters and Schleicher \cite{MorSch25} show that the maximal component $\cluster_{\max}(n)$ satisfies that $\log{|\cluster_{\max}(n)|}/\log{n}\convas \tfrac{1}{2}-\sqrt{(\tfrac{1}{2}-\gamma)^2-\beta(1-2\gamma)}$.

The phase transition results in Theorem \ref{thm-PT-PAM} \cite{HazHofRay23} apply to uniform, as well as preferential, attachment models, and it can be expected that the subcritical results extend to preferential attachment models with major technical adjustments. This is ongoing. A related problem is the so-called {\em Bernoulli} preferential attachment model, studied by Dereich and M\"orters in \cite{DerMor09, DerMor13}. Also in this model, vertices appear dynamically, and, conditionally on the graph $G_n$ at time $n$, vertex $n+1$ connects to vertex $v\in[n]$ {\em independently} with probability $f(\indeg(v,n))/n$, where $\indeg(v,n)$ denotes the in-degree of vertex $v$ at time $n$, i.e., the number of edges connected to $v$ in the time period $v+1, \ldots, n$, and $f\colon \N\to (0,\infty)$ is an appropriate function. Dereich and M\"orters \cite{DerMor13} determine when this model has a giant component  Remarkably, Eckhoff, M\"orters and Ortgiese \cite{EckMorOrt18} show that for super-critical $f$ that are close to critical, the giant is small, in that the proportion of vertices in it is approximately $\e^{-c/\sqrt{\rho-\rho_c}}$, where $\rho_c$ corresponds to the critical edge density, and $\rho>\rho_c$ is a near-critical edge density. The results also apply to near-critical supercritical percolation on such supercritical  graphs. This `small' giant component corresponds to an {\em infinite-order phase transition} that was already found by Riordan \cite{Rior05}, building on \cite{BolJanRio05,BolRio05}, in related settings, including the $m$-out graph, but not yet for preferential attachment models with fixed out-degrees.

\subsection{Proofs for percolation on dynamic random graphs.}
\label{sec-proofs-dynamic}
\myparagraph{The phase transition.}
The proof for the critical percolation threshold in \cite{HazHofRay23} consists of three main steps. First, we show that preferential attachment graphs are so-called {\em large-set expanders}, meaning that even large sets have many edges pointing out of them. This enables us to verify the conditions posed by Alimohammadi, Borgs, and Saberi \cite{AliBorSab23}. Under their conditions, the proportion of vertices in the largest component in a random graph sequence converges to the survival probability of percolation on the local limit. In particular, the critical percolation threshold for both the graph and its local limit are identical.
Second, we identify $1/\pi_c$ as the spectral radius of the mean offspring operator of the P\'olya point tree (PPT), which is the local limit of preferential attachment models.
Lastly, we prove that the critical percolation threshold for the PPT is the inverse of the spectral radius of the mean offspring operator. For positive $\delta$, we use sub-martingales to prove sub-criticality and apply spine decomposition theory to demonstrate super-criticality, completing the third step of the proof. For negative $\delta$, we show that the percolated PPT can be stochastically lower bounded by a supercritical branching process, for all $\pi>0$.

%For $\delta\leq 0$ and any $\pi>0$ instead, we prove that the percolated P\'olya point tree dominates a supercritical branching process, proving that the critical percolation threshold equals $0$.

\myparagraph{Dynamics of components in Theorem \ref{thm-max-comp}.} Again we fix $a=0$ and $m=2$. It will be convenient to use a continuous-time construction $(G^\pi_t)_{t\geq 0}$. Let  $(|\clusterpi{1}{t}|)_{t\geq 0}$ denote the size process of the component containing vertex one. 
For continuous time $t\geq 0$, write the collection of components in $G^{\pi}_t$ as 
	\eqn{
	\label{clusters-UA-discrete}
	\vC(t) = \set{\clusterold \subseteq G^{\pi}_t\colon  \clusterold \text{ component}}.
	}
%With a slight abuse of notation, we denote the set of components at continuous time $t$ by $\vC(t)$, and 
Denote the size of the random graph at time $t$ by $N(t)$. It is not hard to see that $t\mapsto N(t)$ is a {\em Yule process}, i.e., a pure-birth process that grows by 1 at rate $N(t)$. A Yule process satisfies that $N(t)\e^{-t}\convd E$, which is an exponential random variable with parameter 1. Thus, $N(t)$ grows exponentially.

Note that the dynamics of the component $\clusterpi{1}{t}$ are as follows:
\begin{itemize}
    \item[(a)] A new vertex enters the system with two edges and then merges another component $\clusterold\in \vC(t)$ with $\clusterpi{1}{t}$. Thus,  $\clusterold\neq \clusterpi{1}{t} $ merges with $\clusterpi{1}{t}$ leading to the transition $|\clusterpi{1}{t}|\leadsto |\clusterpi{1}{t}| + |\clusterold|+1$ at rate $2\pi^2|\clusterpi{1}{t}||\clusterold|/N(t).$
    \item[(b)] A new vertex enters the system with at least one edge and all its edges connect to $\clusterpi{1}{t}$ leading to the transition $|\clusterpi{1}{t}| \leadsto |\clusterpi{1}{t}| + 1$, which occurs at rate $\pi^2 |\clusterpi{1}{t}|^2/N(t) + 2\pi(1-\pi)|\clusterpi{1}{t}|.$
    \item[(c)] A new vertex enters with no edges connected to it; this leads to no change in the size of $\clusterpi{1}{t}$. 
\end{itemize}
Denote the generator of the associated continuous-time Markov process by $\mathcal{L}$, and, for $k \ge 2$, define the $k^{\rm th}$ susceptibility in continuous time as
	\begin{equation}
	\label{susdef}
    	\csusceptibilitypi{t}{k} :=\frac{1}{N(t)}\sum_{\clusterold \in \vC(t)}|\clusterold|^{k}= \frac{1}{N(t)}\sum_{i=1}^{N(t)}|\clusterpi{i}{t}|^{k-1}, \quad t \ge 0.
	\end{equation} 
Then
	\begin{align}
	\label{dynamics-connected-component-continuous}
    \mathcal{L} |\clusterpi{1}{t}| &= \sum_{\clusterold \in \vC(t)\colon \clusterold\neq \clusterpi{1}{t}}\frac{2\pi^2}{N(t)}|\clusterpi{1}{t}||\clusterold|(|\clusterold|+1) + \left(\pi^2 \frac{|\clusterpi{1}{t}|^2}{N(t)} + 2\pi(1-\pi)|\clusterpi{1}{t}|\right)\\
    &= 2\pi^2|\clusterpi{1}{t}| \left(\csusceptibilitypi{t}{2} - \frac{|\clusterpi{1}{t}|^2}{N(t)}\right) + 2\pi^2|\clusterpi{1}{t}|\left(1 - \frac{|\clusterpi{1}{t}|}{N(t)}\right)+ \left(\pi^2 \frac{|\clusterpi{1}{t}|^2}{N(t)} + 2\pi(1-\pi)|\clusterpi{1}{t}|\right)\nn\\
    &\le \left(2\pi^2\csusceptibilitypi{t}{2} + 2\pi\right)|\clusterpi{1}{t}|.\nn
\end{align}
As a result, for $t\geq 0$, the process 
	\eqn{
	\label{super-martingale}
	M(t)=|\clusterpi{1}{t}| \exp\left(-\int_0^t [2\pi^2 \csusceptibilitypi{u}{2} + 2\pi]~ \dif u \right)
	}
is a non-negative super-martingale, and thus converges almost surely to a finite random variable $M(\infty)$. In order to gain insight into the asymptotics of $|\clusterpi{1}{t}|$, we need to understand the asymptotics of $t\mapsto \csusceptibilitypi{t}{2}$, which we explain next.

\myparagraph{Susceptibilities.}
Define the $k^{\rm th}$ susceptibility functions in discrete time by
\begin{equation}
    \label{eqn:k-suscep}
    \susceptibilitypi{n}{k}= s_k(n)=\frac{1}{n} \sum_{\clusterold\subseteq \vC(n)} |\clusterold|^{k} =  \frac{1}{n} \sum_{v=1}^n |\clusterpi{v}{n}|^{k-1} .
\end{equation}
Note that the discrete and continuous susceptibilities are related through $\susceptibilitypi{N(t)}{k}=\csusceptibilitypi{t}{k}$. From the above, $\susceptibilitypi{n}{2}$ can be viewed as the expected component size of a uniformly chosen vertex in $G^{\pi}_n$ given $G^{\pi}_n$. Although this is a `local' quantity, our proofs will show that this quantity is also key to the asymptotics of `global' quantities like (maximal) component sizes.

\begin{theorem}[Convergence of susceptibility \cite{BanBhaHofRay25}]
\label{thm-suscep-comp}
Fix $a=0$ and $m=2$. For any $\pi < \pi_c$, $\susceptibilitypi{n}{2}\convas \susceptibilitypi{\infty}{2}$, where
	\begin{equation}
    	\susceptibilitypi{\infty}{2} := \frac{(1 - 4\pi)-  \sqrt{8\pi^2 - 8\pi + 1}}{4\pi^2}.
	\label{eqn:s2-def}
	\end{equation}
\end{theorem}
We explain the proof of Theorem \ref{thm-suscep-comp} below. What is remarkable, is that the susceptibility $\susceptibilitypi{\infty}{2}$ in \eqref{eqn:s2-def} remains bounded {\em for all} $\pi<\pi_c$, with a finite limit as $\pi\nearrow \pi_c$. This is unlike related percolation settings, where the expected component size blows up as $\pi\nearrow \pi_c$, and as such is indicative of a different percolation universality class.

\myparagraph{Convergence of components in Theorem \ref{thm-max-comp}: martingales.} Recall \eqref{super-martingale}. Since $\susceptibilitypi{n}{2}\convas \susceptibilitypi{\infty}{2}$, and $\csusceptibilitypi{t}{k}=\susceptibilitypi{N(t)}{k}$ with $N(t)\convas \infty$, it follows that also $\csusceptibilitypi{t}{k}\convas \susceptibilitypi{\infty}{2}$. Thus, almost surely,
	\[
	\int_0^t [2\pi^2 \csusceptibilitypi{u}{2} + 2\pi]~ \dif u
	=t(2\pi^2\susceptibilitypi{\infty}{2}+2\pi)+o(t)=t \alpha(\pi) +o(t),
	\]
by \eqref{eqn:s2-def} and \eqref{eqn:alpha-pi}. As a result, $|\clusterpi{1}{t}|=\e^{t\alpha(\pi)+o(t)}$. Together with $|\clusterpi{1}{t}|$ in continuous time being equal to $|\clusterpi{1}{N(t)}|$ in discrete time, and the fact that $N(t)=\e^{t+o(t)}$, this explain why $|\clusterpi{1}{n}|=n^{\alpha(p)+o(1)}$. To improve this to almost sure convergence of $n\rightarrow |\clusterpi{1}{n}|n^{-\alpha(p)}$, we improve the results to show that $n^{\gamma}(\susceptibilitypi{n}{2}-\susceptibilitypi{\infty}{2})\convas 0$ for some small $\gamma>0$, so that $\int_0^t (\csusceptibilitypi{u}{2}-\susceptibilitypi{\infty}{2})~ \dif u$ converges almost surely as well. To prove that $\zeta_1>0$ almost surely, instead, we consider a martingale $M'(t)=|\clusterpi{1}{t}|^{-1} \exp\left(\int_0^t [2\pi^2 \csusceptibilitypi{u}{2} + 2\pi+R(u)]~ \dif u \right)$, for some appropriate error term $u\mapsto R(u)$. Since the non-negative martingale $M'(t)$ converges almost surely, $\prob_{\pi}(M'(\infty)=\infty)=0$, from which it follows that $\prob_{\pi}(M(\infty)=0)=0$.

\myparagraph{Convergence of susceptibilities in Theorem \ref{thm-suscep-comp}: stochastic approximations.}
%For $\ell \geq 1$, let  $x_\ell^\pi(n)=x_\ell(n)$ denote the number of components of size $\ell$ in $G^{\pi}_n$. For fixed $L\geq 1$, define the $L^{\rm th}$ truncated second susceptibility via 
%	\begin{equation}
%    	\label{eqn:truncated-sus}
%    	\susceptibilitypi{n}{2,L} = \frac{1}{n} \sum_{v=1}^{n}|\clusterpi{v}{n}|\indic{|\clusterpi{v}{n}|\leq L} = \frac{1}{n} \sum_{\ell = 1}^{L} \ell^2 x_{\ell}(n). 
%	\end{equation}
The evolution of the process $(\susceptibilitypi{n}{2})_{n\geq 1}$ in discrete time can be written as 
        \begin{equation}
        \label{SAscheme}
            \susceptibilitypi{n+1}{2} - \susceptibilitypi{n}{2} = \frac{1}{n+1}\left[F(\susceptibilitypi{n}{2}) + R_n + \xi_{n+1} \right], \qquad n\geq 1,
        \end{equation}
where the function $F$ equals $F(s) = 2\pi^2 s^2+ (4\pi-1)s+1$, the error term $R_n$ satisfies $|R_n| \leq K \susceptibilitypi{n}{4}/n$ for a positive constant $K$, and the martingale differences  $(\xi_{n+1})_{n\geq 1}$ given by $\xi_{n+1}=(n+1)\big(\susceptibilitypi{n+1}{2} - \expec[\susceptibilitypi{n+1}{2}\,|\,G^{\pi}_n]\big)$ satisfy $\expec[\xi_{n+1}^2\,|\,G^{\pi}_n]\leq K (\susceptibilitypi{n}{3})^2$. The stochastic recursion in \eqref{SAscheme} follows by an analysis similar to the one in \eqref{dynamics-connected-component-continuous}, applied to the sum of squares of component sizes rather than that of a single vertex. For this, we add vertices one by one, and inspect the evolution of the sum of squares of the components. When two edges are added, they merge components $\clusterold$ and $\clusterold'$ to become a component of size $|\clusterold|+|\clusterold'|+1$ with probability $2|\clusterold||\clusterold'|/n^2$. Summing out over all pairs of components $\clusterold, \clusterold'$ explains why the evolution is {\em quadratic}. The precise shape of $s\mapsto F(s)$ follows by carefully analysing all terms.

To analyse the asymptotics of the stochastic evolution of $\susceptibilitypi{n}{2}$, we note that we can write
	\[
	F(s)=b(s-\lambda_1)(s-\lambda_2),
	\]
where $b=2\pi^2$, $\lambda_1=[1 - 4\pi-  \sqrt{8\pi^2 - 8\pi + 1}]/[4\pi^2]=\susceptibilitypi{\infty}{2}$, and $\lambda_2=[1 - 4\pi+  \sqrt{8\pi^2 - 8\pi + 1}]/[4\pi^2]>\lambda_1.$ The values $\lambda_1$ and $\lambda_2$ correspond to the {\em fixed points} of the dynamics in \eqref{SAscheme}. Since $F'(s)=2bs-b(\lambda_1+\lambda_2)$, we have that $F'(\lambda_1)<0$ and $F'(\lambda_2)>0$, so that $\lambda_1$ is the stable fixed point, while $\lambda_2$ is the unstable one. As such, it comes as no surprise that the dynamics drives $\susceptibilitypi{n}{2}$ to the stable fixed point $\lambda_1=\susceptibilitypi{\infty}{2},$ which explains Theorem \ref{thm-suscep-comp}. The proof follows by a careful analysis of the evolution, using as a starting point that $\susceptibilitypi{n}{2}\convp \susceptibilitypi{\infty}{2}$ by a local convergence argument, combined with a uniform integrability argument in \cite[Chapter 8]{Hofs24}. This argument identifies $\susceptibilitypi{\infty}{2}$ also as the expected component size of the root of percolation on the uniform attachment version of the P\'olya point tree, which is the local limit of our model. Since $\susceptibilitypi{n}{2}\convp \susceptibilitypi{\infty}{2}$, for every $\eta>0$, $n\mapsto \susceptibilitypi{n}{2}$ must enter the $\eta$ neighbourhood of $\susceptibilitypi{\infty}{2}$ infinitely often. Each time, the evolution has a strictly positive probability of never leaving this $\eta$ neighbourhood, so eventually it has to succeed. 

An essential ingredient in the proof is a preliminary estimate that shows that $|\clusterpi{\max}{n}|/\sqrt{n}\convas 0$. Such a bound is crucial, since it shows that the process $n\mapsto \susceptibilitypi{n}{2}$ makes small steps, so that it cannot jump far when it is close to $\susceptibilitypi{\infty}{2}$.

\section{Concluding comments.}
\label{sec-concluding-comments}
In this paper, we discussed recent progress on percolation on random graphs, showing that there are many universality classes. We focussed on the {\em sizes} of the components, thus ignoring their metric structure.  Much progress has been made on the metric structure of critical rank-1 random graphs. See \cite{AddBroGol10, AddBroGol12, BasBhaBroSenWan25, BhaDhaHofSen22, BhaHofSen18, BroDuqWan21, ConGol23, BhaDhaHofSen20, GolHaaSen22} and the references therein. We again refer to \cite[Part III]{Hofs25}, as well as \cite{AddGol24}, for extensive overviews.

Many problems remain, and I expect them to continue attracting attention in the coming years. Of course, the story of near-critical {\em dynamic} random graphs is far from complete, and we do not really understand why they behave so differently compared to the rank-1 settings. So far, {\em spatial} random graphs at criticality have not yet attracted much attention, in part because the method of choice in spatial mean-field-like systems, the lace expansion (see e.g., \cite{HeyHof17} for an extensive introduction), is quite technical, and would need considerable adaptation in order to be able to deal with the significant {\em inhomogeneity} and {\em dependence} present in most spatial random graphs.

\section*{Acknowledgments.} This work is supported in part by the Netherlands Organisation for Scientific Research (NWO) through the Gravitation {\sc Networks} grant 024.002.003.
I thank my collaborators Sayan Banerjee, Shankar Bhamidi, Souvik Dhara, Rajat Hazra, Johan van Leeuwaarden, Rounak Ray, and Sanchayan Sen for the exciting research that forms the basis of this paper, and the joy of discovering it together.
% SIAM recommends using BibTeX
% if using BibTeX
{\small 
\bibliographystyle{siamplain}
%\bibliography{../bib/bibBooks}
%\end{document}
\providecommand{\noopsort}[1]{}\def\cprime{$'$}

}
\end{document}